\documentclass[a4paper,10pt]{article}

\usepackage{fullpage}

\usepackage{natbib}

\usepackage{amsmath}

\usepackage{algorithm2e}
\usepackage{enumerate}
\usepackage{ifthen}
\newboolean{showcomments}
\setboolean{showcomments}{true} 
\ifthenelse{\boolean{showcomments}}
{\newcommand{\nb}[2]{\fbox{\bfseries\small\sffamily#1}{\sf\small$\triangleright$\textbf{#2}$\triangleleft$}}}
{\newcommand{\nb}[2]{}}

\newcommand{\greedy}[1]{\ensuremath{\mathcal G(#1)}}
\newcommand{\tpi}[1]{\ensuremath{T_{#1}}}
\newcommand{\gpi}[1]{\ensuremath{\gamma P_{{#1}}}}
\newcommand{\lgpi}[1]{\ensuremath{\gamma^\ell P_{{#1}}}}
\newcommand{\ident}{\ensuremath{I}}

\newcommand{\norm}[1]{{\ensuremath{\left\lVert #1 \right\rVert}}}
\newcommand{\abs}[1]{{\ensuremath{\left\lvert #1 \right\rvert}}}

\newcommand{\shift}{s}
\newcommand{\bellres}{b}
\newcommand{\distance}{d}
\newcommand{\loss}{l}
\def\T{T}
\def\R{\mathbb{R}}
\def\E{\mathbb{E}}

\usepackage{dsfont}
\newcommand{\indicator}[1]{\ensuremath{\mathds{1}_{[#1]}}}
\newcommand{\gsuml}[1]{\ensuremath{\frac{\gamma^{\ell #1}-\gamma^{\ell(m+1)}}{1-\gamma^\ell}}}

\usepackage{pgf}
\usepackage{tikz}
\usetikzlibrary{arrows,automata,positioning}

\newcommand{\valuefunction}[2]{ 
\begin{tikzpicture}[x=.85cm,y=.85cm,scale=0.95]

\pgfmathsetmacro\states{18}
\pgfmathsetmacro\iter{10}

\pgfmathparse{int(\states/(#1*#2+1)+1)};
\let\n\pgfmathresult

\begin{scope}
\clip (0,0) rectangle (\states,\iter);

\node at (\states-2, \iter-2) {\eqref{eq:1f}};

\ifnum #2>0
  \filldraw[fill=blue!20] (0,0) -- (0, \iter-2/3) -- ++(#1, 0)
     \foreach \i in {1,...,\n} {
       --  ++(1,-1) -- ++(#1 * #2, 0)
     }
     -- (\states,0)
  ;

  \pgfmathparse{\n * #2};
  \let\x\pgfmathresult
  \foreach \i in {0,...,\x} { 
    \node at (\iter/2-.5+#1/2+#1*\i,\iter/2-.5) {\eqref{eq:1d}};
  }
\fi

\filldraw[fill=yellow!20] (0,0) -- (0,\iter-1.5) -- (\iter-1.5,0);
\node at (2, 2) {\eqref{eq:1a}};

\ifnum #2>0
\foreach \i in {0,...,\n} {
  \filldraw[fill=orange!20] (#1 * #2 * \i + \i, \iter-2/3-\i) -- ++(#1, 0) -- ++(0,-2/3) -- ++(-#1,0) -- cycle;
  \node at (#1 /2+#1*#2*\i+\i,\iter-1-\i) {\eqref{eq:1e}};
}
\else
  \filldraw[fill=orange!20] (1/6, \iter-2/3) -- ++(#1-1,0) -- ++(\iter-.5,-\iter+.5) -- ++(1-#1,0) -- cycle;
  \ifnum #1>1
    \node at (3.5+#1 /2,3.5) {\eqref{eq:1e}};
  \else
    \node[rotate=-45] at (\iter/2-.5+#1 /2,\iter/2 -.5) {\eqref{eq:1e}};
  \fi 
\fi

\filldraw[fill=green!20,shift={(#1 -1,0)}] (1/6, \iter-2/3) -- ++(1,0) -- ++(\iter-.5,-\iter+.5) -- ++(-1,0) -- cycle;
\node[rotate=-45] at (\iter/2-.5+#1,\iter/2-.5) {\eqref{eq:1b}};

\filldraw[fill=red!20] (\iter-1.5,0) -- (0,\iter-1.5) -- (0,\iter-2/3) -- (1/6,\iter-2/3) -- (\iter-.5,0);
\node[rotate=-45] at (\iter/2-.5,\iter/2-.5) {\eqref{eq:1c}};

\foreach \ii [evaluate=\ii as \i using \ii-1] in {1,...,\n} {
  \ifnum #2>0
  \foreach \jj [evaluate=\jj as \j using \jj-1] in {1,...,#2} {
    \pgfmathparse{\i * (#1*#2+1) + \j * #1};
    \let\x\pgfmathresult
    \filldraw[fill=green!20,shift={(2*#1,0)}]    
       (\x+ 1/6, \iter-2+1/3-\i) -- ++(1,0) -- ++(\iter-1.5-\i,-\iter+1.5+\i) -- ++(-1,0) -- cycle;
    \node[rotate=-45] at (\iter/2-.5+2*#1+\x-\i, \iter/2-.5) {\eqref{eq:1b}};
  }
  \fi
}

\draw[style=help lines,step=.85cm] (0,0) grid (\states,\iter);
\end{scope}

\pgfmathparse{\states+1}
\let\last\pgfmathresult
\foreach \x in {1,...,\last} \node[anchor=south] at (\x-1,\iter) {\x};
\foreach \y in {0,...,\iter}   \node[anchor=east]  at (0,\iter-\y)     {\y};

\node[above=0.6cm] at (\states/2, \iter) {$i$ (state)};
\node[rotate=90] at (-.8cm,\iter/2) {$k$ (iteration)};

\end{tikzpicture}
}
\usepackage{tikz}
\usetikzlibrary{matrix,decorations.pathreplacing,calc}

\pgfkeys{tikz/mymatrixenv/.style={decoration=brace,every left delimiter/.style={xshift=3pt},every right delimiter/.style={xshift=-3pt}}}
\pgfkeys{tikz/mymatrix/.style={matrix of math nodes,inner sep=2pt,column sep=0.5em,row sep=0.5em,nodes={inner sep=0pt}}}
\pgfkeys{tikz/mymatrixbrace/.style={decorate,thick}}
\newcommand\mymatrixbraceoffseth{0.5em}
\newcommand\mymatrixbraceoffsetv{0.2em}

\newcommand*\mymatrixbraceright[4][m]{
    \draw[mymatrixbrace] ($(#1.north west)!(#1-#3-1.south west)!(#1.south west)-(\mymatrixbraceoffseth,0)$)
        -- node[left=2pt] {#4} 
        ($(#1.north west)!(#1-#2-1.north west)!(#1.south west)-(\mymatrixbraceoffseth,0)$);
}

\newcommand*\mymatrixbracetop[4][m]{
    \draw[mymatrixbrace] ($(#1.north west)!(#1-1-#2.north west)!(#1.north east)+(0,\mymatrixbraceoffsetv)$)
        -- node[above=2pt] {#4} 
        ($(#1.north west)!(#1-1-#3.north east)!(#1.north east)+(0,\mymatrixbraceoffsetv)$);
}

\usepackage[colorlinks=true,breaklinks=true,bookmarks=true,urlcolor=blue,
     citecolor=blue,linkcolor=blue,bookmarksopen=false,draft=false]{hyperref}

\usepackage{amsthm}
\usepackage{amscd}
\usepackage{amsfonts}
\usepackage{authblk}

\newcommand{\appdx}[1]{Section~\ref{#1}}

\def\beginproof{\begin{proof}}
\def\endproofb{\end{proof}}

\newtheorem{theorem}{Theorem}
\newtheorem{lemma}{Lemma}
\newtheorem{definition}{Definition}
\newtheorem{corollary}{Corollary}
\newtheorem{example}{Example}

\title{Tight Performance Bounds for \\ Approximate Modified Policy Iteration \\ with Non-Stationary Policies}

\author[1]{Boris Lesner\thanks{lesnerboris@gmail.com}}
\author[1,2]{Bruno Scherrer\thanks{bruno.scherrer@inria.fr}}
\affil[1]{INRIA Nancy Grand Est, Team MAIA}
\affil[2]{CNRS: UMR7503 – Universit\'e de Lorraine}

\begin{document}

\maketitle

\begin{abstract}%
We consider approximate dynamic programming for the infinite-horizon
stationary $\gamma$-discounted optimal control problem formalized by
Markov Decision Processes. While in the exact case it is known that
there always exists an optimal policy that is stationary, we show that
when using value function approximation, looking for a non-stationary
policy may lead to a better performance guarantee. We define a
non-stationary variant of MPI that unifies a broad family of
approximate DP algorithms of the literature. For this algorithm we
provide an error propagation analysis in the form of a performance
bound of the resulting policies that can improve the usual performance
bound by a factor $O\left(1-\gamma\right)$, which is significant when
the discount factor $\gamma$ is close to 1.  Doing so, our approach
unifies recent results for Value and Policy Iteration. Furthermore, we
show, by constructing a specific deterministic MDP, that our
performance guarantee is tight.

\end{abstract}

\section{Introduction}

We consider a discrete-time dynamic system whose state transition
depends on a control. We assume that there is a {state space}
$X$. When at some state, an action is chosen from a finite {action space}
$A$. The current state $x \in X$ and action $a \in A$ characterizes through a homogeneous probability kernel $P(dx|x,a)$  the next state's distribution. At each transition, the system is given a reward $r(x,a,y) \in \R$ where $r:X \times A \times Y \to \R$ is the instantaneous {reward function}. In this context, we aim at determining a sequence of actions that maximizes the expected discounted sum of rewards from any starting state $x$:
$$
\E \left[ \left. \sum_{k=0}^{\infty}\gamma^k r(x_k,a_k,x_{k+1})  ~\right|~ x_0=x,~ x_{t+1} \sim P(\cdot|x_t,a_t),~ a_0,~ a_1,~ \dots \right],
$$
where $0<\gamma < 1$ is a discount factor. 
The tuple $\langle X, A, P, r, \gamma \rangle$ is called a \emph{Markov Decision Process} (MDP) and the associated optimization problem \emph{infinite-horizon stationary discounted optimal control}~\citep{puterman1994markov,bertsekas1996neuro} .

An important result of this setting is that there exists at least one
stationary deterministic policy, that is a function $\pi:X \rightarrow A$ that maps states into actions,  that is optimal~\citep{puterman1994markov}. As a consequence, the problem is usually recast as looking for the stationary deterministic policy $\pi$ that maximizes for all initial state $x$ the quantity
\begin{equation}
\label{bellvdef}
v_\pi(x):= \E \left[ \left. \sum_{k=0}^{\infty}\gamma^k r(x_k,\pi(x_k),x_{k+1}) ~\right|~ x_0=x,~ x_{t+1} \sim P_\pi(\cdot|x_t)\right],
\end{equation}
also called  the {value of policy $\pi$} at state $x$, and where we wrote $P_\pi(dx|x)$ for the stochastic kernel $P(dx|x,\pi(s))$ that chooses actions according to policy $\pi$. We shall similarly write $r_\pi:X \to \R$ for the function 
that giving the immediate reward while following policy $\pi$:
$$
\forall x,~~r_\pi(x)= \E \left[r(x_0,\pi(x_0),x_1)~|~ x_1 \sim P_\pi(\cdot|x_0) \right].
$$
Two linear operators are associated to the stochastic kernel $P_\pi$: a left operator on functions
$$
\forall f \in \R^X,~~\forall x \in X,~~(P_\pi f)(x) = \int f(y)P_\pi(dy|x) = \E \left[ f(x_1) ~|~ x_1 \sim P_\pi(\cdot|x_0) \right],
$$
and a right operator on distributions:
$$
\forall \mu,~~(\mu P_\pi)(dy) = \int P_\pi(dy|x) \mu(dx).
$$
In words, $P_\pi f(x)$ is the expected value of $f$ after following policy $\pi$ for a single time-step starting from $x$, and $\mu P_\pi$ is the distribution of states after a single time-step starting from $\mu$.

Given a policy $\pi$, it is well known that the value $v_\pi$ is the unique solution of the following Bellman equation:
$$
v_\pi = r_\pi + \gamma P_\pi v_\pi.
$$
In other words, $v_\pi$ is the fixed point of the affine operator $\T_\pi v := r_\pi + \gamma P_\pi v$. 

The {optimal value} starting from state $x$ is defined as
$$
v_*(x):= \max_\pi v_\pi(x).
$$
It is also well known that $v_*$ is characterized by the
following Bellman equation:
$$
v_* = \max_\pi (r_\pi + \gamma P_\pi v_*) = \max_\pi \T_\pi v_*,
$$
where the \mbox{max} operator is componentwise.
In other words, $v_*$ is the fixed point of the nonlinear  operator $\T
v:=\max_{\pi}\T_\pi v$. For any value vector
$v$, we say that a policy $\pi$ is {greedy with respect to the value $v$} if it satisfies:
$$
\pi \in \arg\max_{\pi'} \T^{\pi'} v
$$
or equivalently $\T_\pi v = \T v$. We write, with some abuse of
notation\footnote{There might be several policies that are greedy with respect to some value $v$.} $\greedy{v}$ any policy that is greedy with respect to
$v$. The notions of optimal value function and greedy policies are
fundamental to optimal control because of the following property: any
policy $\pi_*$ that is greedy with respect to the optimal value is an
{optimal policy} and its value $v_{\pi_*}$ is equal to $v_*$.

\par

Given an MDP, we consider approximate versions of the Modified Policy Iteration (MPI) algorithm
\citep{puterman1978modified}.  Starting from an arbitrary value function~$v_0$, MPI generates a sequence of value-policy pairs
\begin{align}
\pi_{k+1} & =  \greedy{v_k} & & \mbox{(greedy step)} \nonumber \\
v_{k+1} & =  (T_{\pi_{k+1}})^{m+1} v_k + \epsilon_k & &\mbox{(evaluation step)} \nonumber
\end{align}
where $m \geq 0$ is a free parameter. At each iteration $k$, the term
$\epsilon_k$ accounts for a possible approximation in the evaluation
step.  MPI generalizes the well-known dynamic programming algorithms
Value Iteration (VI) and Policy Iteration (PI) for values $m=0$ and
$m=\infty$, respectively. In the exact case ($\epsilon_k=0$), MPI requires
less computation per iteration than PI (in a way similar to VI) and
enjoys the faster convergence (in terms of number of iterations) of
PI~\citep{puterman1978modified,puterman1994markov}. 

It was recently shown that controlling the errors $\epsilon_k$ when running MPI is
sufficient to ensure some performance
guarantee~\citep{scherreropi,scherrer2012approximate,scherrer2012,canbolat2012}.
For instance, we have the following performance bound, that is remarkably independent of the
parameter $m$.
\begin{theorem}[{\citet[Remark 2]{scherrer2012approximate}}]
\label{thm:standardbound}
Consider MPI with any parameter $m \ge 0$. Assume there exists an
$\epsilon>0$ such that the errors satisfy
$\|\epsilon_k\|_\infty<\epsilon$ for all $k$. Then, the \emph{loss}
due to running policy $\pi_k$ instead of the optimal policy $\pi_*$
satisfies
\[
\norm{v_* - v_{\pi_k}}_\infty \leq \frac {2(\gamma-\gamma^k)}{(1-\gamma)^2} \epsilon+\frac{2\gamma^k}{1-\gamma}\norm{v_*-v_0}_\infty.
\]
\end{theorem}
In the specific case corresponding to VI ($m=0$) and PI ($m=\infty$),
this bound matches performance guarantees that have been known for a
long time \citep{singh94,bertsekas1996neuro}.  The constant
$\frac{2\gamma}{(1-\gamma)^2}$ can be very big, in particular when
$\gamma$ is close to $1$, and consequently the above bound is commonly
believed to be conservative for practical applications.
Unfortunately, this bound cannot be improved: \citet[Example
  6.4]{bertsekas1996neuro} showed that the bound is tight for PI,
\citet{ScherrerLesner2012} proved that it is tight for VI\footnote{Though the MDP instance used to show the tightness of the bound for VI is the same as that for PI~\citep[Example 6.4]{bertsekas1996neuro}, \citet{ScherrerLesner2012} seem to be the first to argue about it in the literature.}, and we will
prove in this article\footnote{Theorem~\ref{thm:tight} page~\pageref{thm:tight} with $\ell=1$.} the---to our knowledge unknown---fact that it
is also tight for MPI. In other words, improving the performance bound
requires to change the algorithms.

\section{Main Results}

Even though the theory of optimal control states that there exists a
stationary policy that is optimal, \citet{ScherrerLesner2012} recently showed
that the performance bound of Theorem~\ref{thm:standardbound} could be improved in the specific cases 
$m=0$ and $m=\infty$ by considering
variations of VI  and PI that build \emph{periodic non-stationary
  policies} (instead of stationary policies). In this article, we
consider an original MPI algorithm that generalizes these variations
of VI and PI (in the same way the standard MPI algorithm generalizes
standard VI and PI). Given some free parameters $m \ge 0$ and $\ell
\ge 1$, an arbitrary value function $v_0$ and an arbitrary set of
$\ell-1$ policies $\pi_{0}, \pi_{-1}, \pi_{-\ell+2}$, consider the algorithm
that builds a sequence of value-policy pairs as follows:
\begin{align}
\pi_{k+1} & =  \greedy{v_k} & & \mbox{(greedy step)} \nonumber \\
v_{k+1} & =  (T_{\pi_{k+1,\ell}})^{m} T_{\pi_{k+1}} v_k + \epsilon_k. & &\mbox{(evaluation step)} \nonumber
\end{align}
While the greedy step is identical to the one of the standard MPI algorithm, the evaluation step involves two new objects that we describe now. $\pi_{k+1,\ell}$ denotes the 
  {periodic non-stationary policy that loops in reverse order on the last $\ell$ generated policies}:
\[
\pi_{k+1,\ell} = \underbrace{\pi_{k+1}\ \pi_{k}\ \cdots\ \pi_{k-\ell+2}}_{{\tiny \mbox{last $\ell$ policies}}}~\underbrace{\pi_{k+1}\ \pi_{k}\ \cdots\ \pi_{k-\ell+2}}_{{\tiny \mbox{last $\ell$ policies}}}\cdots
\]
Following the policy $\pi_{k+1,\ell}$ means that the first action is selected by $\pi_{k+1}$, the second one
by $\pi_{k}$, until the $\ell^\mathrm{th}$ one by $\pi_{k-\ell+2}$,
then the policy loops and the next actions are selected by $\pi_{k+1}$, $\pi_{k}$, so on and so forth.
In the above algorithm,  $T_{\pi_{k+1,\ell}}$ is the linear
Bellman operator associated to $\pi_{k+1,\ell}$:
\[
T_{\pi_{k+1,l}} = T_{\pi_{k+1}} T_{\pi_k} \dots T_{\pi_{k-\ell+2}},
\]
that is the operator of which the unique fixed point is the value function of $\pi_{k+1,\ell}$.
After $k$ iterations, the output of the algorithm is the periodic non-stationary policy $\pi_{k,\ell}$. 

For the values $m=0$ and $m=\infty$, one respectively recovers the
variations of VI\footnote{As already noted by \cite{ScherrerLesner2012}, the only difference between this variation of VI and the standard VI algorithm is what is output by the algorithm. Both algorithms use the very same evaluation step: $v_{k+1}=T_{\pi_{k+1}}v_k$. However, after $k$ iterations, while standard VI returns the last stationary policy $\pi_k$, the variation of VI returns the non-stationary policy $\pi_{k,\ell}$.} and PI recently proposed by
\cite{ScherrerLesner2012}. When $\ell=1$, one recovers the standard
MPI algorithm by \cite{puterman1978modified} (that itself generalizes the standard VI and PI algorithm).
As it generalizes all previously proposed algorithms, we will simply refer to this new algorithm as MPI with parameters $m$ and $\ell$.

On the one hand, using this new algorithm
may require more memory since one must store $\ell$ policies instead of one. On the other hand, our first main result, proved in \appdx{proof:bound}, shows that this extra memory allows to improve the performance guarantee.
\begin{theorem}
  \label{thm:bound}
  Consider MPI with any parameters $m \ge 0$ and $\ell \ge 1$.
  Assume there exists an
$\epsilon>0$ such that the errors satisfy
$\|\epsilon_k\|_\infty<\epsilon$ for all $k$. Then, the \emph{loss}
due to running policy $\pi_{k,\ell}$ instead of the optimal policy $\pi_*$
satisfies
  \[
  \norm{v_*-v_{\pi_{k,\ell}}}_\infty \leq
  \frac{2(\gamma-\gamma^{k})}{(1-\gamma)(1-\gamma^\ell)}\epsilon +
  \frac{2\gamma^k}{1-\gamma}\norm{v_*-v_0}_\infty.
\]
\end{theorem}
As already observed for the standard MPI algorithm, this performance
bound is independent of $m$. For any $\ell \ge 1$, it is a factor
$\frac{1-\gamma}{1-\gamma^\ell}$ better than in
Theorem~\ref{thm:standardbound}. 
Using $\ell=\left\lceil \frac 1{1-\gamma}\right\rceil$ yields\footnote{
Using the facts that $1-\gamma \leq -\log \gamma$ and $\log\gamma\leq 0$, we have
$
\log\gamma^{\ell} \leq \log\gamma^{\frac 1 {1-\gamma}} \leq \frac 1 {-\log\gamma}\log\gamma = -1
$
hence $\gamma^{\ell} \leq e^{-1}$. Therefore
$
 \frac 2{1-\gamma^\ell} \leq \frac 2 {1-e^{-1}} < 3.164.
$
}
 a performance bound of
\[
  \norm{v_*-v_{\pi_{k,\ell}}}_\infty <
  \frac{3.164(\gamma-\gamma^{k})}{1-\gamma}\epsilon +
  \frac{2\gamma^k}{1-\gamma}\norm{v_*-v_0}_\infty,
  \]
and constitutes asymptotically an improvement of order $O(1-\gamma)$, which is significant when $\gamma$ is close to 1.
In fact, Theorem~\ref{thm:bound} is a
generalization of Theorem~\ref{thm:standardbound} for $\ell > 1$ (the
bounds match when $\ell=1$).  While this result was obtained
through two independent proofs for the variations of VI and PI proposed
by~\cite{ScherrerLesner2012}, the more general setting that we
consider here involves a unified proof that extends that provided for the
standard MPI ($\ell=1$) by \cite{scherrer2012}. Moreover, our result is much more general since it applies to all the variations of MPI for any $\ell$ and $m$.

\begin{figure}[ht]
\begin{center}
\resizebox{\linewidth}{!}{
\begin{tikzpicture}[->,>=stealth',shorten >=1pt,auto,node distance=2.1cm,
                    semithick,every loop/.style={<-,shorten <=1pt}]
\tikzstyle{every state}=[thick,minimum size=1.2cm]
\tikzstyle{cont}=[label={center:$\dots$}, auto, minimum size=1cm]
\node[state] (S1)  {$1$};
\node[state] (S2)[right of=S1] {$2$};
\node[state] (S3)[right of=S2] {$3$};

\node[cont]  (S_)[right of=S3] {};

\node[state] (SL1) [right of=S_] {$\ell$};
\node[state] (SL2) [right of=SL1] {$\ell+1$};
\node[state] (SL3) [right of=SL2] {$\ell+2$};

\node[cont]  (S__) [right of=SL3] {};

\path
(S1) edge [loop above] node {$0$} (S1)
(S2) edge [bend left] node {$-2\gamma\epsilon$} (SL2)
(S3) edge [bend left] node {$-2(\gamma+\gamma^2)\epsilon$} (SL3)
;

\path
(S2)  edge node {$0$} (S1)
(S3)  edge node {$0$} (S2)
(S_)  edge node {$0$} (S3)
(SL1) edge node {$0$} (S_)
(SL2) edge node {$0$} (SL1)
(SL3) edge node {$0$} (SL2)
(S__) edge node {$0$} (SL3)
;
\end{tikzpicture}}
\end{center}
\caption{The deterministic MDP matching the bound of Theorem~\ref{thm:bound}.}
\label{fig:mdp}
\end{figure}
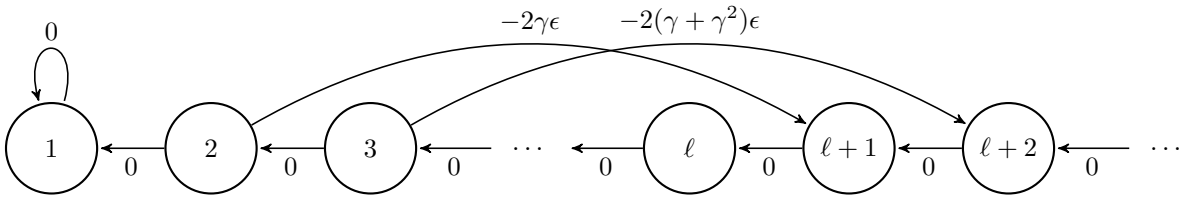
The second main result of this article, proved in \appdx{sec:tight}, is that the bound of Theorem~\ref{thm:bound} is tight, in the precise sense formalized by the following theorem.
\begin{theorem}
\label{thm:tight}
For all parameter values $m \ge 0$ and $\ell \ge 1$, for all $\epsilon > 0$, there exists an MDP instance, an initial value function
$v_0$, a set of initial policies $\pi_0,\pi_{-1},\dots,\pi_{-\ell+2}$ and a sequence of error terms $(\epsilon_k)_{k \ge 1}$ satisfying $\|\epsilon_k\|_\infty \le \epsilon$, such
that for all iterations $k$, the bound of
Theorem~\ref{thm:bound} is satisfied with equality.
\end{theorem}
This theorem generalizes the (separate) tightness results for
PI~\citep{bertsekas1996neuro} and for VI~\citep{ScherrerLesner2012}
where the problem constructed to attain the bound is a specialization
of the one we use in \appdx{sec:tight}. To our knowledge, this result
is new even for the standard MPI algorithm ($m$ arbitrary but $\ell=1$),
and for the non-trivial non-stationary variations of VI
($m=0$, $\ell>1$) and PI ($m=\infty$, $\ell>1$).  The proof considers
a generalization of the MDP instance used to prove the tightness of the bound for
VI~\citep{ScherrerLesner2012} and PI~\citep[Example 6.4]{bertsekas1996neuro}.
Precisely, this MDP consists of states $\{1,2,\dots \}$, two actions:
left ($\leftarrow$) and right ($\rightarrow$); the reward function $r$ and transition kernel $P$ are characterized as
follows for any state $i \ge 2$:
\begin{align*}
r(i,\leftarrow) &= 0,     & r(i, \rightarrow) &= -2\frac{\gamma-\gamma^i}{1-\gamma}\epsilon, \\
P(i|i+1,\leftarrow) &= 1, & P(i+\ell-1|i,\rightarrow) &= 1,
\end{align*}
and $r(1)=0$ and $P(1|1)1$ for state $1$ (all the other transitions having zero probability mass).
As a shortcut, we will use the notation $r_i$ for the non-zero reward $r(i,\rightarrow)$ in state $i$.
Figure~\ref{fig:mdp} depicts the general structure of this MDP.  It is
easily seen that the optimal policy $\pi_*$ is to take $\leftarrow$ in
all states $i \ge 2$, as doing otherwise would incur a negative reward. Therefore,
the optimal value $v_*(i)$ is $0$ in all states $i$. 
The proof of the above theorem considers that we run MPI  with $v_0=v_*=0$, 
$\pi_0=\pi_{-1}=\dots=\pi_{\ell+2}=\pi_*$, 
and the following sequence of error terms:
\[
\forall i,~~
\epsilon_k(i) = \left\{\begin{array}{ll}
-\epsilon&\text{if } i=k,\\
\epsilon&\text{if } i=k+\ell,\\
0&\text{otherwise}.
\end{array}\right.
\]
In such a case, one can prove that the sequence of 
policies $\pi_1,\pi_2,\dots,\pi_k$ that are generated up to iteration $k$ is such
that for all $i \le k$, the policy $\pi_i$ takes $\leftarrow$ in all states but $i$, where it takes $\rightarrow$.  As a
consequence, a non-stationary policy $\pi_{k,\ell}$ built from this
sequence takes $\rightarrow$ in $k$ (as dictated by $\pi_k$), which
transfers the system into state $k+\ell-1$ incurring a reward of
$r_k$. Then the policies $\pi_{k-1}, \pi_{k-2}, \dots, \pi_{k-\ell+1}$
are followed, each indicating to take $\leftarrow$ with $0$
reward. After $\ell$ steps, the system is again is state $k$ and, by the
periodicity of the policy, must again use the action $\pi_k(k) =\
\rightarrow$. The system is thus stuck in a loop, where every $\ell$
steps a negative reward of $r_k$ is received. Consequently, the value of 
this policy from state $k$ is:
\[
v_{\pi_{k,\ell}}(k) = r_k + \gamma^\ell (r_k + \gamma^\ell
(r_k+\cdots)\cdots) = \frac{r_k}{1-\gamma^\ell} =
-\frac{\gamma-\gamma^{k}}{(1-\gamma)(1-\gamma^\ell)}2\epsilon.
\] 
As a consequence, we get the following lower bound, 
\[
\norm{v_* - v_{\pi_{k,\ell}}}_\infty \ge |v_{\pi_{k,\ell}}(k)| =
\frac{\gamma-\gamma^{k}}{(1-\gamma)(1-\gamma^\ell)}2\epsilon
\]
which \emph{exactly} matches the upper bound of
Theorem~\ref{thm:bound} (since $v_0=v_*=0$).  The proof of this result
involves computing the values $v_k(i)$ for all states $i$, steps $k$
of the algorithm, and values $m$ and $\ell$ of the parameters, and
proving that the policies $\pi_{k+1}$ that are greedy with respect to
these values satisfy what we have described above.  Because of the
cyclic nature of the MDP, the shape of the value function is quite
complex---see for instance Figures~\ref{fig:value-l2-m3} and~\ref{fig:value-l3-m2} in
\appdx{sec:tight}---and the exact derivation is tedious.  For clarity,
this proof is deferred to \appdx{sec:tight}.

\section{Discussion}

Since it is well known that there exists an optimal policy that is
stationary, our result---as well as those of
\citet{ScherrerLesner2012}---suggesting to consider non-stationary
policies may appear strange. There exists, however, a very simple
approximation scheme of discounted infinite-horizon control
problems---that has to our knowledge never been documented in the
literature---that sheds some light on the deep reason why
non-stationary policies may be an interesting option.  Given an
infinite-horizon problem, consider approximating it by a
finite-horizon discounted control problem by ``cutting the horizon''
after some sufficiently big instant $T$ (that is assume there is no
reward after time $T$). Contrary to the original infinite-horizon
problem, the resulting finite-horizon problem is non-stationary, and
has therefore \emph{naturally} a non-stationary solution that is built
by dynamic programming in reverse order. Moreover, it can be shown
\citep[by adapting the proof of Theorem 2.5.1]{kakade} that solving
this finite-horizon with VI with a potential error of $\epsilon$ at
each iteration, will induce at most a performance error of $2
\sum_{i=0}^{T-1}
\gamma^t\epsilon=\frac{2(1-\gamma^T)}{1-\gamma}\epsilon$. If we add
the error due to truncating the horizon ($\gamma^T
\frac{\max_{s,a}|r(s,a)|}{1-\gamma}$), we get an overall error of
order $O\left(\frac{1}{1-\gamma}\epsilon\right)$ for a memory $T$ of
the order of\footnote{ We use the fact that $\gamma^T K
  <\frac{\epsilon}{1-\gamma} ~\Leftrightarrow~T > \frac{\log
    \frac{(1-\gamma)K}{\epsilon}}{\log \frac{1}{\gamma}} \simeq
  \frac{\log \frac{(1-\gamma)K}{\epsilon}}{1-\gamma}$ with
  $K=\frac{\max_{s,a}|r(s,a)|}{1-\gamma}$.}  $\tilde
O\left(\frac{1}{1-\gamma}\right)$.  Though this approximation scheme
may require a significant amount of memory (when $\gamma$ is close to
$1$), it achieves the same $O(1-\gamma)$ improvement over the standard
MPI performance bound as our MPI new scheme proposed through our
generalization of MPI with two parameters $m$ and $\ell$.  In
comparison, the new proposed algorithm can be seen as a more flexible
way to make the trade-off between the memory and the quality.

A practical limitation of Theorem~\ref{thm:bound} is that it assumes
that the errors $\epsilon_k$ are controlled in max norm. In practice,
the evaluation step of dynamic programming algorithm is usually done
through some regression scheme---see for instance
\citep{bertsekas1996neuro,antos2007fitted,antos2007value,scherrer2012approximate}---and
thus controlled through some weighted quadratic $L_{2,\mu}$ norm,
defined as $\|f\|_{2,\mu}=\sqrt{\int f(x) \mu(dx)}$.
\citet{munos2003,munos2007performance} originally developed such
analyzes for VI and PI. \citet{FaMuSz10} and
\citet{scherrer2012approximate} later improved it. Using a technical
lemma due to \citet[Lemma 3]{scherrer2012approximate}, one can easily
deduce\footnote{Precisely, Lemma 3 of \citep{scherrer2012approximate}
  should be applied to Equation~\eqref{eq:abs-loss}
  page~\pageref{eq:abs-loss} in \appdx{proof:bound}.} from our
analysis (developed in \appdx{proof:bound}) the following performance
bound.
\begin{corollary}
Consider MPI with any parameters $m \ge 0$ and $\ell \ge 1$.
  Assume there exists an
$\epsilon>0$ such that the errors satisfy
$\|\epsilon_k\|_{2,\mu}<\epsilon$ for all $k$. Then, the expected (with respect to some initial measure $\rho$) \emph{loss}
due to running policy $\pi_{k,\ell}$ instead of the optimal policy $\pi_*$
satisfies 
  \[
  \E\left[v_*(s)-v_{\pi_{k,\ell}}(s)~|~ s \sim \rho \right] \leq
  \frac{2(\gamma-\gamma^{k})C_{1,k,\ell}}{(1-\gamma)(1-\gamma^\ell)} \epsilon +
  \frac{2\gamma^kC_{k,k+1,1}}{1-\gamma}\norm{v_*-v_0}_{2,\mu},
\]
where 
$$
C_{j,k,l}=\frac{(1-\gamma)(1-\gamma^l)}{\gamma^j-\gamma^k}\sum_{i=j}^{k-1} \sum_{n=i}^\infty \gamma^{i + l n}  c(i + l n)
$$
is a convex combination of concentrability coefficients based on Radon-Nikodym derivatives
$$
c(j)= \max_{\pi_1,\cdots,\pi_j}\left\| \frac{d (\rho P_{\pi_1}P_{\pi_2}\cdots P_{\pi_j})}{d \mu} \right\|_{2,\mu}.
$$
\end{corollary}
With respect to the previous bound in max norm, this bound involves
extra constants $C_{j,k,l} \ge 1$. Each such coefficient $C_{j,k,l}$
is a convex combination of terms $c(i)$, that each quantifies the
difference between 1)~the distribution $\mu$ used to control the
errors and 2)~the distribution obtained by starting from $\rho$ and
making $k$ steps with arbitrary sequences of policies. Overall, this
extra constant can be seen as a measure of stochastic smoothness of
the MDP (the smoother, the smaller). Further details on these
coefficients can be found in
\citep{munos2003,munos2007performance,FaMuSz10}.

The next two sections contain the proofs of our two main results, that are Theorem~\ref{thm:bound} and~\ref{thm:tight}.

\section{Proof of Theorem~\ref{thm:bound}}
\label{proof:bound}

Throughout this proof we will write $P_k$ (resp. $P_*$) for the transition kernel $P_{\pi_k}$ (resp. $P_{\pi_*}$) induced by the stationary policy $\pi_k$ (resp. $\pi_*$). We will write $T_k$ (resp. $\tpi{*}$) for the associated Bellman operator. Similarly, we will write $P_{k,\ell}$ for the transition kernel associated with the non-stationary policy $\pi_{k,\ell}$ and $\tpi{k,\ell}$ for its associated Bellman operator.

For $k\geq 0$ we define the following quantities:
\begin{itemize}
\item $\bellres_k = \tpi{k+1}v_{k} - \tpi{k+1,\ell}\tpi{k+1}v_{k}$. This
  quantity which we will call the \emph{residual} may be viewed as a
  non-stationary analogue of the Bellman residual $v_k-\tpi{k+1}v_{k}$.
\item $\shift_k = v_k - v_{\pi_{k,\ell}} - \epsilon_k$. We will call it
  \emph{shift}, as it measures the shift between the value $v_{\pi_{k,\ell}}$ and the estimate $v_k$ before incurring the error.
\item $\distance_k = v_* - v_k + \epsilon_k$. This quantity, called
    \emph{distance} thereafter, provides the distance between the $k^\text{th}$ value function (before the error is added)  and the optimal value function.
  \item $\loss_k = v_* - v_{\pi_{k,\ell}}$. This is the \emph{loss} of
    the policy $v_{\pi_{k,\ell}}$. The loss is always non-negative since
    no policy can have a value greater than or equal to $v_*$.
\end{itemize}
The proof is outlined as follows. We first provide a bound on
$\bellres_k$ which will be used to express both the bounds on $\shift_k$
and $\distance_k$. Then, observing that $\loss_k = \shift_k+\distance_k$
will allow to express the bound of $\norm{\loss_k}_\infty$ stated by
Theorem~\ref{thm:bound}. Our arguments extend those made by 
\citet{scherrer2012approximate} in the specific case $\ell=1$.

We will repeatedly use the fact that since policy $\pi_{k+1}$ is greedy with respect to $v_k$, we have
\begin{equation}
\label{eq:g}
\forall \pi',~~ T_{k+1} v_k \ge T_{\pi'} v_k.
\end{equation}
For a non-stationary policy $\pi_{k,\ell},$ the
induced $\ell$-step transition kernel is
\[
P_{k,\ell} = P_kP_{k-1}\cdots P_{k-\ell+1}.
\]
As a consequence, for any function $f:\mathcal S\rightarrow\mathbb{R}$, the operator $\tpi{k,\ell}$ may be expressed as:
\begin{equation*}
\tpi{k,\ell}f =r_k+\gamma P_{k,1}r_{k-1}+\gamma^2 P_{k,2}r_{k-2} + \cdots +\gamma^{\ell-1} P_{k,\ell-1}r_{k-\ell+1} + \gamma^{\ell}P_{k,\ell}f
\end{equation*}
then, for any function $g:\mathcal S\rightarrow\mathbb{R}$,  we have
\begin{equation}
\label{eq:f-g}
\tpi{k,\ell}f - \tpi{k,\ell}g = \lgpi{k,\ell}(f-g)
\end{equation}
and
\begin{equation}
\label{eq:f+g}
  \tpi{k,\ell}(f+g) = \tpi{k,\ell}f + \lgpi{k,\ell}(g).
\end{equation}

The following notation will be useful.
\begin{definition}[\citet{scherrer2012approximate}]
  For a positive integer $n$, we define $\mathbb{P}_n$ as the set of discounted
  transition kernels that are defined as follows:
  \begin{enumerate}
  \item for any set of $n$ policies $\{\pi_1,\dots,\pi_n\}$, $(\gamma
    P_{\pi_1})(\gamma P_{\pi_2})\cdots(\gamma P_{\pi_n})
    \in\mathbb{P}_n$,
  \item for any $\alpha \in (0,1)$ and $P_1, P_2 \in \mathbb{P}_n$,
    $\alpha P_1 + (1-\alpha)P_2 \in \mathbb{P}_n$
  \end{enumerate}
With some abuse of notation, we write $\Gamma^n$ for denoting any element of $\mathbb{P}_n$.
\end{definition}
\begin{example}[$\Gamma^n$ notation]
  If we write a transition kernel $P$ as
  $P=\alpha_1\Gamma^i+\alpha_2\Gamma^j\Gamma^k =
  \alpha_1\Gamma^i+\alpha_2\Gamma^{j+k}$, it should be read as: ``There
  exists $P_1\in\mathbb{P}_i$,$P_2\in\mathbb{P}_j$,$P_3\in\mathbb{P}_k$
  and $P_4\in\mathbb{P}_{j+k}$ such that $P=\alpha_1P_1+\alpha_2P_2P_3 =
  \alpha_1P_1+\alpha_2P_4$.''.
\end{example}

We first provide three lemmas bounding the residual, the shift and the
distance, respectively.

\begin{lemma}[residual bound]
\label{lem:bellres}
The residual $\bellres_k$ satisfies the following bound:
\[
\bellres_k \leq \sum_{i=1}^{k}\Gamma^{(\ell m+1)(k-i)}x_{i} + \Gamma^{(\ell m+1)k}\bellres_0
\]
where 
\[
x_k = (\ident -\Gamma^\ell)\Gamma\epsilon_k.
\]
\end{lemma}
\beginproof
We have:
\begin{align*}
\bellres_{k} &= \tpi{k+1}v_{k} - \tpi{k+1,\ell}\tpi{k+1}v_{k}\\
&\leq \tpi{k+1}v_k - \tpi{k+1,\ell}\tpi{k-\ell+1}v_k&\hspace{-2cm}\{\tpi{k+1}v_k\geq \tpi{k-\ell+1}v_k~\eqref{eq:g} \}\\
&=\tpi{k+1}v_k - \tpi{k+1}\tpi{k,\ell}v_k\\
&=\gpi{k+1}\left( v_k - \tpi{k,\ell}v_k\right)\\
&=\gpi{k+1}\left( (\tpi{k,\ell})^m \tpi k v_{k-1}+\epsilon_k - \tpi{k,\ell}\left((\tpi{k,\ell})^m \tpi k v_{k-1}+\epsilon_k\right)\right)\\
&=\gpi{k+1}\left( (\tpi{k,\ell})^m \tpi k v_{k-1} -(\tpi{k,\ell})^{m+1} \tpi k v_{k-1} + (\ident-\lgpi {k,\ell})\epsilon_k\right)&\{\eqref{eq:f+g}\}\\
&=\gpi{k+1}\left( (\lgpi {k,\ell})^m \left(\tpi k v_{k-1} -\tpi{k,\ell} \tpi k v_{k-1}\right) + (\ident-\lgpi {k,\ell})\epsilon_k\right)&\{\eqref{eq:f-g}\}\\
&=\gpi{k+1}\left( (\lgpi {k,\ell})^m \bellres_{k-1} + (\ident-\lgpi {k,\ell})\epsilon_k\right).\\
\end{align*}
Which can be written as 
\[
\bellres_k \leq \Gamma(\Gamma^{\ell m}\bellres_{k-1} + (I-\Gamma^\ell)\epsilon_k) = \Gamma^{\ell m+1}\bellres_{k-1} + x_k.
\]
Then, by induction:
\[
\bellres_k \leq \sum_{i=0}^{k-1}\Gamma^{(\ell m+1)i}x_{k-i} + \Gamma^{(\ell m+1)k}\bellres_0 = \sum_{i=1}^{k}\Gamma^{(\ell m+1)(k-i)}x_{i} + \Gamma^{(\ell m+1)k}\bellres_0.
\]
\endproofb
\begin{lemma}[distance bound]
\label{lem:distance}
The distance $\distance_k$ satisfies the following bound:
\[
\distance_k\leq\sum_{i=1}^{k}\sum_{j=0}^{mi-1}\Gamma^{\ell j +i-1}x_{k-i} +\sum_{i=1}^k\Gamma^{i-1}y_{k-i}+ z_k,
\]
where
\[
y_k = -\Gamma\epsilon_k
\]
and
\[
z_k = \sum_{i=0}^{mk-1} \Gamma^{k-1+\ell i}\bellres_0 + \Gamma^k\distance_0.
\]
\end{lemma}
\beginproof
First expand $\distance_k$:
\begin{align*}
\distance_k &= v_* - v_k + \epsilon_k\\
&= v_* - (\tpi{k,\ell})^m \tpi k v_{k-1}\\
\begin{split}
  &= v_* - \tpi k v_{k-1} + \tpi k v_{k-1} - \tpi{k,\ell} \tpi k v_{k-1} +\tpi{k,\ell} \tpi k v_{k-1}- (\tpi{k,\ell})^2\tpi k v_{k-1}\\[.3cm]
  &\quad  + (\tpi{k,\ell})^2 \tpi k v_{k-1} - \cdots - (\tpi{k,\ell})^{m-1} \tpi k v_{k-1} + (\tpi{k,\ell})^{m-1} \tpi k v_{k-1} - (\tpi{k,\ell})^m \tpi k v_{k-1}\hspace{-1cm}
\end{split}\\
&= v_* - \tpi k v_{k-1} + \sum_{i=0}^{m-1} (\tpi{k,\ell})^i\tpi k v_{k-1} - (\tpi{k,\ell})^{i+1}\tpi k v_{k-1}\\
&= \tpi * v_* - \tpi k v_{k-1} + \sum_{i=0}^{m-1} (\lgpi {k,\ell})^i\left(\tpi k v_{k-1} - \tpi{k,\ell} \tpi k v_{k-1}\right)
&\{\eqref{eq:f-g}\}\\
&\leq \tpi * v_* - \tpi * v_{k-1} + \sum_{i=0}^{m-1} (\lgpi {k,\ell})^i\bellres_{k-1}&\{\tpi k v_{k-1} \geq \tpi * v_{k-1}~\eqref{eq:g}\}\\
&=\gpi *(v_* - v_{k-1}) + \sum_{i=0}^{m-1} (\lgpi {k,\ell})^i\bellres_{k-1}
&\{\eqref{eq:f-g}\}\\
&=\gpi *\distance_{k-1} - \gpi *\epsilon_{k-1}+ \sum_{i=0}^{m-1} (\lgpi {k,\ell})^i\bellres_{k-1}
&\{\distance_k=v_* - v_k + \epsilon_k\}\\
&=\Gamma\distance_{k-1} + y_{k-1} + \sum_{i=0}^{m-1} \Gamma^{\ell i}\bellres_{k-1}.\\
\end{align*}
Then, by induction
\[
\distance_k\leq\sum_{j=0}^{k-1}\Gamma^{k-1-j}\left(y_j+\sum_{p=0}^{m-1} \Gamma^{\ell p}\bellres_{j}\right) + \Gamma^k\distance_0.
\]
Using the bound on $\bellres_k$ from Lemma~\ref{lem:bellres} we get:
\begin{align*}
\distance_k&\leq \sum_{j=0}^{k-1}\Gamma^{k-1-j}\left(y_j+\sum_{p=0}^{m-1} \Gamma^{\ell p}\left(\sum_{i=1}^{j}\Gamma^{(\ell m+1)(j-i)}x_{i} + \Gamma^{(\ell m+1)j}\bellres_0\right)\right) + \Gamma^k\distance_0\\
&=\sum_{j=0}^{k-1}\sum_{p=0}^{m-1} \sum_{i=1}^{j}\Gamma^{k-1-j+\ell p+(\ell m+1)(j-i)}x_{i} + \sum_{j=0}^{k-1}\sum_{p=0}^{m-1}\Gamma^{k-1-j+\ell p+(\ell m+1)j}\bellres_0 + \Gamma^k\distance_0+\sum_{i=1}^k\Gamma^{i-1}y_{k-i}.
\end{align*}
First we have:
\begin{align*}
\sum_{j=0}^{k-1}\sum_{p=0}^{m-1} \sum_{i=1}^{j}\Gamma^{k-1-j+\ell p+(\ell m+1)(j-i)}x_{i}&=\sum_{i=1}^{k-1}\sum_{j=i}^{k-1}\sum_{p=0}^{m-1}\Gamma^{k-1+\ell(p+mj)-i(\ell m +1)}x_{i}\\
&=\sum_{i=1}^{k-1}\sum_{j=0}^{m(k-i)-1}\Gamma^{ k-1 +\ell(j+mi) -i(\ell m +1)}x_{i}\\
&=\sum_{i=1}^{k-1}\sum_{j=0}^{m(k-i)-1}\Gamma^{\ell j + k-i-1}x_{i}\\
&=\sum_{i=1}^{k-1}\sum_{j=0}^{mi-1}\Gamma^{\ell j +i-1}x_{k-i}.\\
\end{align*}
Second we have:
\[
\sum_{j=0}^{k-1}\sum_{p=0}^{m-1}\Gamma^{k-1-j+\ell p+(\ell m+1)j}\bellres_0 = \sum_{j=0}^{k-1}\sum_{p=0}^{m-1}\Gamma^{k-1+\ell (p+ mj)}\bellres_0 = \sum_{i=0}^{mk-1} \Gamma^{k-1+\ell i}\bellres_0 = z_k- \Gamma^k\distance_0.
\]
Hence
\[
\distance_k\leq\sum_{i=1}^{k}\sum_{j=0}^{mi-1}\Gamma^{\ell j +i-1}x_{k-i} +\sum_{i=1}^k\Gamma^{i-1}y_{k-i}+ z_k.
\]
\endproofb
\begin{lemma}[shift bound]
\label{lem:shift}
The shift $\shift_k$ is bounded by:
\[
\shift_k \leq \sum_{i=1}^{k-1}\sum_{j=mi}^\infty\Gamma^{\ell j + i-1}x_{k-i} + w_k,
\]
where
\[
w_k = \sum_{j=mk}^\infty\Gamma^{\ell j+k-1}\bellres_0.
\]
\end{lemma}
\beginproof
Expanding $\shift_k$ we obtain:
\begin{align*}
\shift_k &= v_k - v_{\pi_{k,\ell}} - \epsilon_k\\
&= (\tpi{k,\ell})^m\tpi k v_{k-1} -  v_{\pi_{k,\ell}}\\
&= (\tpi{k,\ell})^m\tpi k v_{k-1} - (\tpi{k,\ell})^\infty\tpi{k,\ell}\tpi k v_{k-1}&\{\forall f:\ v_{\pi_{k,\ell}} = (\tpi{k,\ell})^\infty f\}\\
&= (\lgpi{k,\ell})^m\sum_{j=0}^\infty(\lgpi {k,\ell})^j\left(\tpi k v_{k-1} - \tpi{k,\ell}\tpi k v_{k-1}\right)\\
&= \Gamma^{\ell m} \sum_{j=0}^\infty\Gamma^{\ell j} \bellres_{k-1}\\
&= \sum_{j=0}^\infty\Gamma^{\ell m+\ell j}\bellres_{k-1}.
\end{align*}
Plugging the bound on $\bellres_k$ of Lemma~\ref{lem:bellres} we get:
\begin{align*}
  \shift_k &\leq \sum_{j=0}^\infty\Gamma^{\ell m+\ell j}\left(\sum_{i=1}^{k-1}\Gamma^{(\ell m+1)(k-1-i)}x_{i} + \Gamma^{(\ell m+1)(k-1)}\bellres_0\right)\\
  &= \sum_{j=0}^\infty\sum_{i=1}^{k-1}\Gamma^{\ell m+\ell j + (\ell m+1)(k-1-i)}x_{i} +  \sum_{j=0}^\infty\Gamma^{\ell m+\ell j+(\ell m+1)(k-1)}\bellres_0\\
  &= \sum_{j=0}^\infty\sum_{i=1}^{k-1}\Gamma^{\ell(j+mi)+i-1}x_{k-i} +  \sum_{j=0}^\infty\Gamma^{\ell (j+mk)+k-1}\bellres_0\\
  &= \sum_{i=1}^{k-1}\sum_{j=mi}^\infty\Gamma^{\ell j+i-1}x_{k-i} +  \sum_{j=mk}^\infty\Gamma^{\ell j+k-1}\bellres_0\\
  &= \sum_{i=1}^{k-1}\sum_{j=mi}^\infty\Gamma^{\ell j + i-1}x_{k-i} + w_k.
\end{align*}
\endproofb
\begin{lemma}[loss bound]
\label{lem:loss}
The loss $\loss_k$ is bounded by: 
\[
\loss_k \leq \sum_{i=1}^{k-1}\Gamma^{i}\left(\sum_{j=0}^\infty\Gamma^{\ell j}(\ident-\Gamma^\ell)-\ident\right)\epsilon_{k-i} +\eta_k,
\]
where
\[
\eta_k = z_k + w_k = \sum_{i=0}^{mk-1} \Gamma^{k-1+\ell i}\bellres_0 + \Gamma^k\distance_0 + \sum_{j=mk}^\infty\Gamma^{\ell j+k-1}\bellres_0 = \sum_{i=0}^{\infty} \Gamma^{\ell i+k-1}\bellres_0 + \Gamma^k\distance_0.
\]
\end{lemma}
\beginproof
Using Lemmas~\ref{lem:distance} and~\ref{lem:shift}, we have:
\begin{align*}
\loss_k &= \shift_k + \distance_k\\
&\leq\sum_{i=1}^{k-1}\sum_{j=mi}^\infty\Gamma^{\ell j + i-1}x_{k-i}  + \sum_{i=1}^{k-1}\sum_{j=0}^{mi-1}\Gamma^{\ell j +i-1}x_{k-i}+\sum_{i=1}^k\Gamma^{i-1}y_{k-i}+ z_k +w_k\\
&=\sum_{i=1}^{k-1}\sum_{j=0}^\infty\Gamma^{\ell j + i-1}x_{k-i} + \sum_{i=1}^k\Gamma^{i-1}y_{k-i} + \eta_k.
\end{align*}
Plugging back the values of $x_k$ and $y_k$ and using the fact that $\epsilon_0=0$ we obtain:
\begin{align*}
\loss_k &\leq \sum_{i=1}^{k-1}\sum_{j=0}^\infty\Gamma^{\ell j + i-1}(\ident-\Gamma^\ell)\Gamma\epsilon_{k-i} + \sum_{i=1}^{k-1}\Gamma^{i-1}(-\Gamma)\epsilon_{k-i} - \Gamma^k\epsilon_0+\eta_k\\
&= \sum_{i=1}^{k-1}\left(\sum_{j=0}^\infty\Gamma^{\ell j + i}(\ident-\Gamma^\ell)\epsilon_{k-i} -\Gamma^{i}\epsilon_{k-i}  \right)+\eta_k\\
&= \sum_{i=1}^{k-1}\Gamma^{i}\left(\sum_{j=0}^\infty\Gamma^{\ell j}(\ident-\Gamma^\ell)-\ident\right)\epsilon_{k-i} +\eta_k.\\
\end{align*}
\endproofb
We now provide a bound of $\eta_k$ in terms of $\distance_0$:
\begin{lemma}
\label{lem:eta}
\[
  \eta_k \leq \Gamma^k\left(\sum_{i=0}^{\infty}\Gamma^i(\Gamma - \ident)+\ident\right)\distance_0.
\]
\end{lemma}
\beginproof
First recall that
\[
\eta_k = \sum_{i=0}^{\infty} \Gamma^{\ell i+k-1}\bellres_0 + \Gamma^k\distance_0.
\]
In order to bound $\eta_k$ in terms of $\distance_0$ only, we express $\bellres_0$ in terms of  $\distance_0$:
\begin{align*}
  \bellres_0 &= \tpi 1 v_0 - (\tpi 1)^\ell \tpi 1 v_0\\
  &= \tpi 1 v_0 - (\tpi 1)^2v_0 + (\tpi 1)^2v_0 - \cdots - (\tpi 1)^\ell v_0 +  (\tpi 1)^\ell v_0  - (\tpi 1)^{\ell+1}v_0\\
  &=\sum_{i=1}^\ell(\gpi 1)^i(v_0-\tpi 1 v_0)\\
  &= \sum_{i=1}^\ell(\gpi 1)^i(v_0 - v_* +\tpi * v_* - \tpi * v_0 + \tpi * v_0 -\tpi 1 v_0)\\
  &\leq  \sum_{i=1}^\ell(\gpi 1)^i(v_0 - v_* +\tpi * v_* - \tpi * v_0) & \{ T_1 v_0 \ge T_* v_0 ~ \eqref{eq:g} \}\\
  &= \sum_{i=1}^\ell(\gpi 1)^i(\gpi * - \ident)d_0.
\end{align*}

Consequently, we have:
\begin{align*}
  \eta_k &\leq \sum_{i=0}^{\infty} \Gamma^{\ell i+k-1}\sum_{j=1}^\ell(\gpi 1)^j(\gpi * - \ident)d_0 + \Gamma^k\distance_0\\
  &= \sum_{i=0}^{\infty} \Gamma^{\ell i+k}\sum_{j=0}^{\ell-1}(\gpi 1)^j(\gpi * - \ident)d_0 + \Gamma^k\distance_0\\
  &=\Gamma^k\left(\sum_{i=0}^{\infty} \Gamma^{\ell i}\sum_{j=0}^{\ell-1}\Gamma^j(\Gamma - \ident)+\ident\right)\distance_0\\
  &=\Gamma^k\left(\sum_{i=0}^{\infty}\Gamma^i(\Gamma - \ident)+\ident\right)\distance_0.\\
\end{align*}
\endproofb
We now conclude the proof of Theorem~\ref{thm:bound}.
Taking the absolute value in Lemma~\ref{lem:eta} we obtain:
\[
\abs{\eta_k} \leq \Gamma^k\left(\sum_{i=0}^{\infty}\Gamma^i(\Gamma +
  \ident)+\ident\right)\abs{\distance_0} =
2\sum_{i=k}^{\infty}\Gamma^i\abs{\distance_0}
\]
Since $\loss_k$ is non-negative, from Lemma~\ref{lem:loss} we have:
\begin{equation}
\label{eq:abs-loss}
  \abs{\loss_k} \leq
  \sum_{i=1}^{k-1}\Gamma^{i}\left(\sum_{j=0}^\infty\Gamma^{\ell
      j}(\ident+\Gamma^\ell)+\ident\right)\abs{\epsilon_{k-i}}
  +\abs{\eta_k}
=2\sum_{i=1}^{k-1}\Gamma^{i}\sum_{j=0}^\infty\Gamma^{\ell
      j}\abs{\epsilon_{k-i}}
  +2\sum_{i=k}^{\infty}\Gamma^i\abs{\distance_0}.
\end{equation}
Since $\norm{v}_\infty = \max \abs{v}$, $\distance_0 = v_* - v_0$ and
$\loss_k = v_* - v_{\pi_{k,\ell}}$, we can take the maximum in \eqref{eq:abs-loss} and conclude that:
\[
\norm{v_* - v_{\pi_{k,\ell}}}_\infty \leq \frac{2(\gamma-\gamma^k)}{(1-\gamma)(1-\gamma^\ell)}2\epsilon + \frac{\gamma^k}{1-\gamma}\norm{v_*-v_0}_\infty.
\]

\section{Proof of Theorem~\ref{thm:tight}}
\label{sec:tight}

We shall prove the following result.
\begin{lemma}
\label{lem:policyvalue}
Consider MPI with parameters $m \ge 0$ and $\ell \ge 1$ applied on the problem of Figure~\ref{fig:mdp}, starting from $v_0=0$ and all initial policies $\pi_0,\pi_{-1},\dots,\pi_{-\ell+2}$ equal to $\pi_*$. Assume that at each iteration $k$, the following error terms are applied, for some $\epsilon\geq 0$:
\[
\forall i,~~
\epsilon_k(i) = \left\{\begin{array}{ll}
-\epsilon&\text{if } i=k\\
\epsilon&\text{if } i=k+\ell\\
0&\text{otherwise}
\end{array}\right..
\]
Then MPI can\footnote{We write here ``can'' since at each iteration, several policies will be greedy with respect to the current value.} generate a sequence of value-policy pairs that is described below.

For all iterations $k\geq 1$, the policy $\pi_k$ always takes the optimal action in all states, that is
\begin{equation}
  \label{eq:policy}
\forall i \ge 2,~~
  \pi_k(i) = \left\{\begin{array}{ll}
      \rightarrow &\text{\rm if } i = k\\
      \leftarrow&\text{\rm otherwise}
    \end{array}\right.
\end{equation}

\addtocounter{equation}{1}
For all iterations $k\geq 1$, the value function $v_k$ is defined as follows:
\begin{itemize}
\item For all $i<k$:
\begin{equation}
 \label{eq:1a}
v_k(i) = -\gamma^{(k-1)(\ell m+1)}\epsilon
\tag{\arabic{equation}.a}
\end{equation}
\item For all $i$ such that $k\leq i \leq k+((k-1)m+1)\ell$:
\begin{itemize}
  \item For $i = k+(qm+p+1)\ell$ with $q\geq 0$ and $0\leq p < m$ (i.e. $i=k+n\ell$, $n\geq 1$):
    \begin{equation}
      \label{eq:1b}
    v_k(i) = \gamma^{q(\ell m+1)}\left( \gsuml{(p+1)}r_{k-q}
      + \indicator{p=0}\epsilon
      + \sum_{j=1}^{k-q-1}\gamma^{i(\ell m+1)}\left(\gsuml{} r_{k-q-j}+\epsilon  \right) \right)
    \tag{\arabic{equation}.b}
    \end{equation}
  \item For $i=k$:
    \begin{equation}
      \label{eq:1c}
    v_k(k) = v_k(k+\ell) + r_k -2\epsilon
    \tag{\arabic{equation}.c}
    \end{equation}
  \item For $i=k+q\ell +p$ with $0\leq q \leq (k-1)m-1$ and $1\leq p < \ell$:
    \begin{equation}
      \label{eq:1d}
      v_k(i) = -\gamma^{(k-1)(\ell m+1)}\epsilon
      \tag{\arabic{equation}.d}
    \end{equation}
  \item Otherwise, i.e. when $i=k+(k-1)m\ell+p$ with  $1\leq p < \ell$:
    \begin{equation}
      \label{eq:1e}
    v_k(i) = 0
    \tag{\arabic{equation}.e}
    \end{equation}
\end{itemize}
\item For all $i > k+((k-1)m+1)\ell$
\begin{equation}
 \label{eq:1f}
  v_k(i) = 0
  \tag{\arabic{equation}.f}
\end{equation}
\end{itemize}

\end{lemma}

\begin{figure}
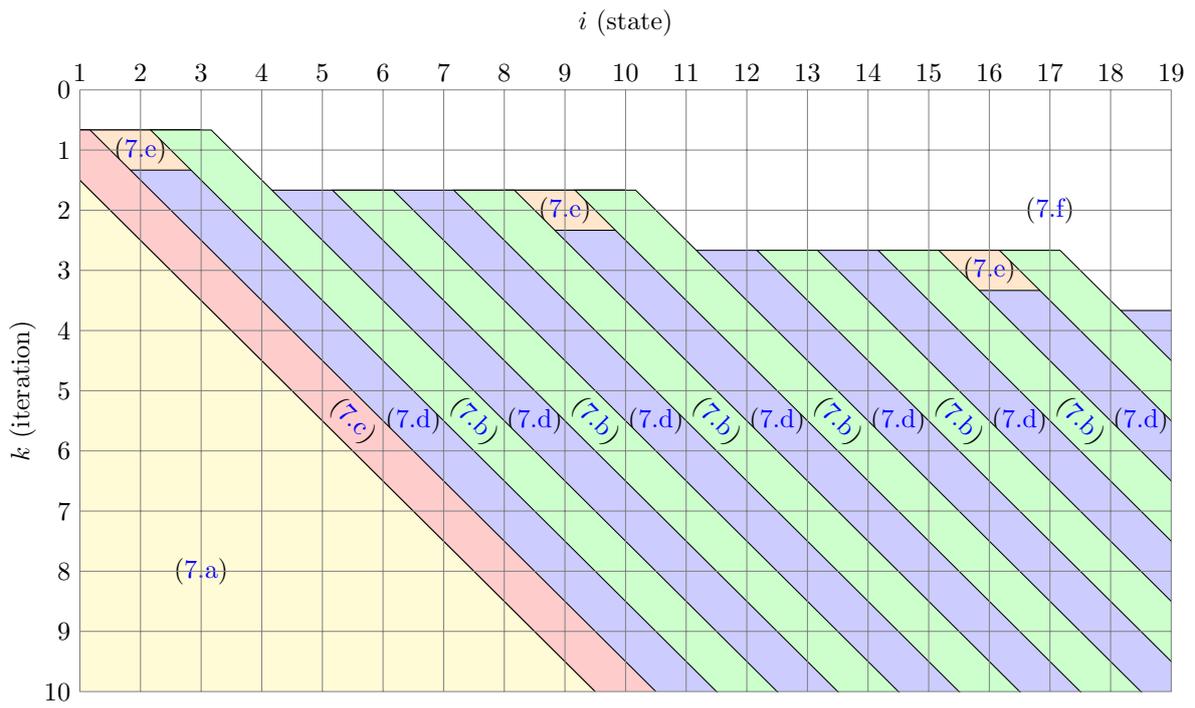

\resizebox{\linewidth}{!}{\valuefunction{2}{3}}
\caption{Shape of the value function with $\ell=2$ and $m=3$.}
\label{fig:value-l2-m3}
\end{figure}

\begin{figure}
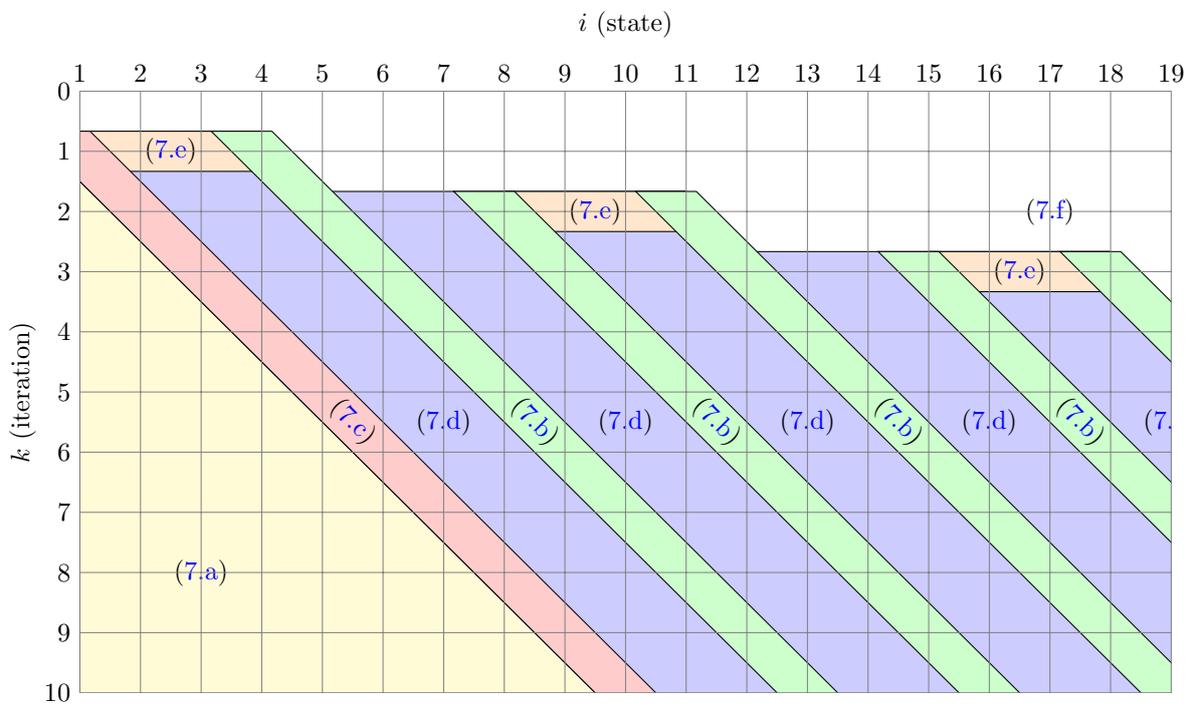

\resizebox{\linewidth}{!}{\valuefunction{3}{2}}
\caption{Shape of the value function with $\ell=3$ and $m=2$.}
\label{fig:value-l3-m2}
\end{figure}

The relative complexity of the different expressions of $v_k$ in
Lemma~\ref{lem:policyvalue} is due to the presence of nested periodic patterns
in the shape of the value fonction along the state space and the
horizon. Figures~\ref{fig:value-l2-m3} and~\ref{fig:value-l3-m2} give
the shape of the value function for different values of $\ell$ and $m$,
exhibiting the periodic patterns.
The proof of Lemma~\ref{lem:policyvalue} is done by recurrence on $k$.

\subsection{Base case $k=1$ }
Since $v_0 = 0$, $\pi_1$ is the optimal policy that takes $\leftarrow$
in all states as desired. Hence, $(\tpi{1,\ell})^m\tpi 1 v_0 = 0$ in all
states. Accounting for the errors $\epsilon_1$ we have $v_1 =
(\tpi{1,\ell})^m\tpi 1 v_0+ \epsilon_1 = \epsilon_1$. As can be seen on
Figures~\ref{fig:value-l2-m3} and~\ref{fig:value-l3-m2}, when $k=1$ we
only need to consider equations \eqref{eq:1b}, \eqref{eq:1c},
\eqref{eq:1e} and \eqref{eq:1f} since the others apply to an empty set
of states.
\par
First, we have 
\[
v_1(1+\ell) = \epsilon_1(1+\ell) = \epsilon
\] which is \eqref{eq:1b} when $q=(k-1)=0$ and $p=0$.  
\par
Second, we have
\[
v_1(1) = \epsilon_1(1) = -\epsilon = \epsilon + 0 - 2\epsilon = v_1(1+\ell) + r_1 - 2\epsilon
\]
which corresponds to \eqref{eq:1c}.
\par
Third, for $1\leq p <\ell$ we have
\[
v_1(1+p) = \epsilon_1(1+p) = 0
\]
corresponding to \eqref{eq:1e}.
\par
Finally, for all the remaining states $i>1+\ell$, we have
\[
v_1(i) = \epsilon_1(i) = 0
\]
corresponding to \eqref{eq:1f}.

The base case is now proved.

\subsection{Induction Step}
We assume that Lemma~\ref{lem:policyvalue}
holds for some \emph{fixed} $k\geq 1$, we now show that it also holds
for $k+1$. 

\subsubsection{The policy $\pi_{k+1}$}
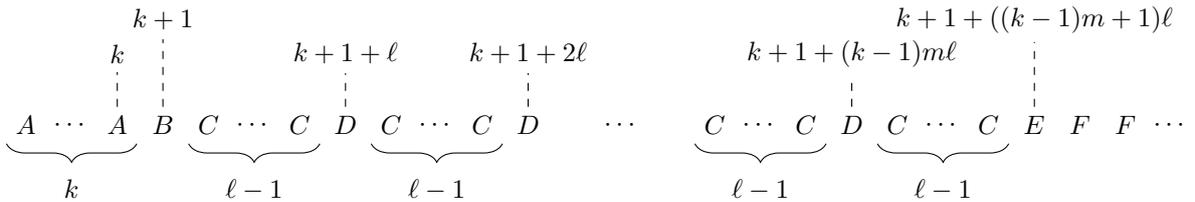
\begin{figure}[ht]
\resizebox{\linewidth}{!}{
\begin{tikzpicture}[node distance=0.6cm]
\usetikzlibrary{decorations.pathreplacing}
\node (S1)                 {$A$};
\node (S1...)[right of=S1] {$\cdots$};
\node (Sk)   [right of=S1...] {$A$};
\node (Sk1)  [right of=Sk] {$B$};

\node (Sk_up) [above of=Sk,yshift=.35cm] {$k$};
\node (Sk1_up) [above of=Sk1,yshift=.8cm] {$k+1$};

\draw[style=dashed] (Sk) -- (Sk_up); 
\draw[style=dashed] (Sk1) -- (Sk1_up); 

\node (Sk2)  [right of=Sk1] {$C$};
\node (Sk2...)  [right of=Sk2] {$\cdots$};
\node (Sk3)  [right of=Sk2...] {$C$};
\node (Sk4)  [right of=Sk3] {$D$};

\node (Sk4_up) [above of=Sk4,yshift=.35cm] {$k+1+\ell$};
\draw[style=dashed] (Sk4) -- (Sk4_up); 

\node (Sk5)  [right of=Sk4] {$C$};
\node (Sk5...)  [right of=Sk5] {$\cdots$};
\node (Sk6)  [right of=Sk5...] {$C$};
\node (Sk7)  [right of=Sk6] {$D$};

\node (Sk7_up) [above of=Sk7,yshift=.35cm] {$k+1+2\ell$};
\draw[style=dashed] (Sk7) -- (Sk7_up); 

\node (Sk7c)  [right=0.6cm of Sk7] {$\cdots$};

\node (Sk8)  [right=0.6cm of Sk7c] {$C$};
\node (Sk8...)  [right of=Sk8] {$\cdots$};
\node (Sk9)  [right of=Sk8...] {$C$};
\node (Sk10)  [right of=Sk9] {$D$};

\node (Sk10_up) [above of=Sk10,yshift=.35cm] {$k+1+(k-1)m\ell$};

\draw[style=dashed] (Sk10) -- (Sk10_up); 

\node (Sk11)  [right of=Sk10] {$C$};
\node (Sk11...)  [right of=Sk11] {$\cdots$};
\node (Sk12)  [right of=Sk11...] {$C$};
\node (Sk13)  [right of=Sk12] {$E$};

\node (Sk13_up) [above of=Sk13,yshift=.8cm] {$k+1+((k-1)m+1)\ell$};
\draw[style=dashed] (Sk13) -- (Sk13_up); 

\node (Sk14)  [right of=Sk13] {$F$};
\node (Sk15)  [right of=Sk14] {$F$};
\node (Sk15...)  [right of=Sk15] {$\cdots$};

\draw [decorate,decoration={brace,amplitude=7pt}]
(Sk.south east) -- (S1.south west) node [black,midway,below,yshift=-10pt] {$k$};

\draw [decorate,decoration={brace,amplitude=7pt}]
(Sk3.south east) -- (Sk2.south west) node [black,midway,below,yshift=-10pt] {$\ell-1$};
\draw [decorate,decoration={brace,amplitude=7pt}]
(Sk6.south east) -- (Sk5.south west) node [black,midway,below,yshift=-10pt] {$\ell-1$};
\draw [decorate,decoration={brace,amplitude=7pt}]
(Sk9.south east) -- (Sk8.south west) node [black,midway,below,yshift=-10pt] {$\ell-1$};

\draw [decorate,decoration={brace,amplitude=7pt}]
(Sk12.south east) -- (Sk11.south west) node [black,midway,below,yshift=-10pt] {$\ell-1$};
\end{tikzpicture}
}
\caption{Policy cases, each state is represented by a letter
  corresponding to a case of the policy $\pi_{k+1}$. Starting from $1$,
  state number increase from left to right.}
\label{fig:policy-cases}
\end{figure}

We begin by showing that the policy $\pi_{k+1}$ is greedy with respect to $v_k$.
Since there is no choice in state $1$ is $\rightarrow$,  we turn our attention to the
other states. There are many cases to consider, each one of them
corresponding to one or more states. These cases, labelled from A through F, are summarized as
follows, depending on the state $i$:
\begin{enumerate}[\qquad(A)]
  \item $1 <i < k+1$
  \item $i=k+1$
  \item $i = k+1+q\ell +p$ with $1\leq p <\ell$ and $0\leq q \leq (k-1)m$
  \item $i = k+1+(qm+p+1)\ell$ with $0\leq p < m$ and $0\leq q < k-1$
  \item $i = k+1+((k-1)m+1)\ell$
  \item $i > k+1+((k-1)m+1)\ell$
  \end{enumerate}
Figure~\ref{fig:policy-cases} depicts how those cases cover the whole state space.
\par
For all states $i>1$ in each of the above cases, we consider the
\emph{action-value functions} $q^\rightarrow_{k+1}(i)$
(resp. $q^\leftarrow_{k+1}(i)$) of action $\rightarrow$
(resp. $\leftarrow$) defined as:
\[
  q^\rightarrow_{k+1}(i) = r_i + \gamma v_k(i-1)\qquad\text{and}\qquad
  q^\leftarrow_{k+1}(i)= \gamma v_k(i+\ell-1).
\]
In case $i=k+1$ (B) we will show that
$q^\rightarrow_{k+1}(i)=q^\leftarrow_{k+1}(i)$ meaning that a policy
$\pi_{k+1}$ greedy for $v_k$ may be either $\pi_{k+1}(k+1)=\
\rightarrow$ or $\pi_{k+1}(k+1)=\ \leftarrow$. In all other cases we
show that $q^\rightarrow_{k+1}(i)<q^\leftarrow_{k+1}(i)$ which implies
that for those $i\neq k+1$, $\pi_{k+1}(i)=\ \leftarrow$, as required by Lemma~\ref{lem:policyvalue}.

\paragraph{A: In states $1<i<k+1$}
We have $q^\rightarrow_{k+1}(i) = r_i+\gamma v_k(i+\ell-1)$ and
$q^\leftarrow_{k+1}(i) = \gamma v_k(i-1)$, depending on the value of $i+\ell-1$, which is
reached by taking the $\rightarrow$ action, we need to consider two
cases:
\begin{itemize}
\item Case 1: $i+\ell-1\neq k$. In this case $v_k(i+\ell-1)$ is
  described by either \eqref{eq:1a} or \eqref{eq:1d} when $i+\ell-1$ is
  less than, or greater than $k$, respectively. In either case we have
$v_k(i+\ell-1)=  -\gamma^{(k-1)(\ell m+1)}\epsilon = v_k(i-1)$ and hence:
    \[
    q^\rightarrow_{k+1}(i) = r_i + \gamma v_k(i+\ell-1) = r_i +\gamma v_k(i-1) < \gamma v_k(i-1) = q^\leftarrow_{k+1}(i)
    \]
    which gives $\pi_{k+1}(i) =\ \leftarrow$ as desired.
  \item Case 2: $i+\ell-1=k$. 
    \begin{align*}
      q^\rightarrow_{k+1}(i) &= r_i + \gamma v_k(k) = r_i + \gamma \left(v_k(k+\ell) + r_k - 2\epsilon\right)&\{\eqref{eq:1c}\}\\
      &=\gamma\left(\sum_{j=0}^{k-1}\gamma^{j(\ell m+1)}\left(\gsuml{}r_{k-j}+\epsilon\right) + r_k -2\epsilon\right)&\{\eqref{eq:1b}\}\\
      &\leq \gamma\left(\sum_{j=0}^{k-1}\gamma^{j(\ell m+1)}\epsilon+  r_k -2\epsilon\right)&\{r_{k-j} \leq 0\}\\
      &= \gamma\left(\sum_{j=1}^{k-1}\left(\gamma^{j(\ell m+1)}\epsilon-2\gamma^j\epsilon\right)-\epsilon\right)&\{r_k =  -2\sum_{j=1}^{k-1}\gamma^{j}\epsilon\}\\
      &< -\gamma\epsilon &\{\gamma^{j(\ell m+1)}\epsilon-2\gamma^j\epsilon<0\}\\
      &< \gamma v_k(i-1)&\!\!\!\!\!\!\!\!\!\!\!\!\{v_k(i-1)=-\gamma^{(k-1)(\ell m+1)}\epsilon\ \eqref{eq:1a}\}\\
      &= q^\leftarrow_{k+1}(i)
    \end{align*}
    giving $\pi_{k+1}(i) =\ \leftarrow$ as desired.
\end{itemize}
\paragraph{B: In state $k+1$}
Looking at the action value function $q_{k+1}^\leftarrow$ in state $k+1$,
we observe that:
\begin{align*}
  q_{k+1}^\leftarrow (k+1) &= \gamma v_k(k) = \gamma\left(r_k - 2\epsilon+v_k(k+\ell)\right)&\text{\{\eqref{eq:1c}\}}\\
  &= \gamma r_k -2\gamma\epsilon +\gamma v_k(k+\ell)\\
  &= r_{k+1} + \gamma v_k(k+\ell)&\text{\{$r_{i+1}=\gamma r_i -2\gamma\epsilon$\}}  \\
  & = q_{k+1}^\rightarrow (k+1)
\end{align*}
This means that the algorithm can take $\pi_{k+1}(k+1) =\ \rightarrow$
so as to satisfy Lemma~\ref{lem:policyvalue}. 
\paragraph{C: In states $i=k+1+q\ell+p$}
We restrict ourselves to the cases when $1\leq p< \ell$ and $0\leq q
\leq (k-1)m$. Three cases for the value of $q$ need to be considered:
\begin{itemize}
\item Case 1: $0\leq q < (k-1)m-1$.
We have:
\begin{align*}
  q_{k+1}^\rightarrow (i) &= r_i + \gamma v_k(k+(q+1)\ell+p)\\
  &=r_i + \gamma v_k(k+q\ell+p)&\text{\{\eqref{eq:1d} independent of $q$\}}\\
  &< \gamma v_k(k+q\ell+p)&\{r_i< 0\}\\
  &=q_{k+1}^\leftarrow (i).
\end{align*}
\item Case 2: $q = (k-1)m-1$
\begin{align*}
q_{k+1}^\rightarrow (i) &= r_i + \gamma  v_k(k+(q+1)\ell+p)\\
&=r_i + \gamma 0&\{\eqref{eq:1e}\}\\
&= -2\epsilon\frac{\gamma-\gamma ^{k+1+q\ell+p}}{1-\gamma}\\
&=-2\epsilon\left(\frac{\gamma-\gamma ^{k+q\ell+p}}{1-\gamma} + \gamma^{k+q\ell+p}\right)\\
&< -\gamma^{k+q\ell+p} \epsilon\\
&= -\gamma^{k+(k-1)\ell m-\ell+p}\epsilon &\{q = (k-1)m-1\}\\
&<- \gamma^{k+(k-1)\ell m}\epsilon=- \gamma^{(k-1)(\ell m+1)+1}\epsilon&\{p-\ell < 0\}\\
&= \gamma v_k(k+q\ell+p)&\{\eqref{eq:1d}\}\\
&= q_{k+1}^\leftarrow (i).
\end{align*}
\item Case 3: $q = (k-1)m$
\begin{align*}
q_{k+1}^\rightarrow (i) &= r_i + \gamma  v_k(k+((k-1)m+1)\ell+p)\\
&= r_i + \gamma 0&\{\eqref{eq:1f}\}\\
&= r_i +\gamma v_k(k+((k-1)m)\ell+p)&\{\eqref{eq:1e}\}\\
& = r_i +\gamma v_k(i-1)\\
&< q_{k+1}^\leftarrow (i).&\{r_i < 0\}
\end{align*}
\end{itemize}

\paragraph{D: In states $i = k+1+(qm+p+1)\ell$}
In these states, we have:
\begin{align}
  q_{k+1}^\leftarrow (i) &= \gamma v_k(k+(qm+p+1)\ell)\notag\\[.4cm]
  q_{k+1}^\rightarrow (i) &= r_{i}+\gamma v_k(k+1+(qm+p+1)\ell+\ell-1)\notag\\
  &=r_{i}+\gamma v_k(k+(qm+p+2)\ell).\label{eq:qval}
\end{align}
As for the right-hand side of \eqref{eq:qval} we need to consider two cases:
\begin{itemize}
  \item Case 1: $p+1< m$:
\end{itemize}
In the following, define
\[
x_{k,q} =  \sum_{j=1}^{k-q-1}\gamma^{j(\ell m+1)}\left(\gsuml{} r_{k-q-j}+\epsilon  \right).
\]
Then,
\begin{align}
  \label{eq:case1}
  \lefteqn{q_{k+1}^\rightarrow (i) = r_{i}+\gamma v_k(k+(qm+(p+1)+1)\ell)}\notag\\
  &=r_{i}+\gamma\gamma^{q(\ell m+1)}\left( \gsuml{(p+2)}r_{k-q}
    + \sum_{j=1}^{k-q-1}\gamma^{j(\ell m+1)}\left(\gsuml{} r_{k-q-j}+\epsilon  \right) \right)&\{\eqref{eq:1b}\}\notag\\
&=r_{i}+\gamma^{q(\ell m+1)+1}\left( \left(\gsuml{(p+1)}-\gamma^{\ell(p +1)}\right)r_{k-q}
+x_{k,q}
\right)\notag\\
  &=r_{i}-\gamma^{(qm+p+1)\ell+q+1} r_{k-q}
  + \gamma^{q(\ell m+1)+1}\left( \gsuml{(p+1)}r_{k-q}
+x_{k,q}
\right)\notag\\
  &=r_{i}-\gamma^{i-k+q} r_{k-q} +  \gamma v_k(k+(qm+p+1)\ell)  -\indicator{p=0}\gamma^{q(\ell m+1)+1}\epsilon&\{\eqref{eq:1b}\}\notag\\
  &\leq r_{i}-\gamma^{i-k+q} r_{k-q} + \gamma v_k(k+(qm+p+1)\ell)\notag\\
  &= r_{i}-\gamma^{i-k+q} r_{k-q} +q_{k+1}^\leftarrow (i).
\end{align}
Now, observe that
\begin{align*}
  \gamma^{i-k+q} r_{k-q} &=  - 2\gamma^{i-k+q} \frac{\gamma-\gamma^{k-q}}{1-\gamma}\epsilon\\
  &= - 2\frac{\gamma^{i-k+q+1}-\gamma^{i}}{1-\gamma}\epsilon\\
  &= - 2\frac{\gamma -\gamma +\gamma^{i-k+q+1}-\gamma^{i}}{1-\gamma}\epsilon\\
  &= - 2\frac{\gamma-\gamma^i}{1-\gamma}\epsilon - 2\frac{-\gamma + \gamma^{i-k+q+1}}{1-\gamma} \epsilon\\
  &= r_{i} - r_{i-k+q+1}.
\end{align*}
Plugging this back into \eqref{eq:case1}, we get:
\begin{align*}
q_{k+1}^\rightarrow (i) &\leq r_{i} -  r_{i} +r_{i-k+q+1} +  q_{k+1}^\leftarrow(i)\\
&< q_{k+1}^\leftarrow (i).&\{r_{i-k+q+1}< 0\}
\end{align*}
\begin{itemize}
  \item Case 2: $p+1= m$:
\end{itemize}
Using the fact that $p+1=m$ implies $\gsuml{(p+1)}= \gamma^{\ell m}$ we have:
\begin{align*}
  \lefteqn{q_{k+1}^\rightarrow (i) = r_{i}+\gamma v_k(k+((q+1)m+1)\ell)}\\
   &=r_{i}+\gamma\gamma^{(q+1)(\ell m+1)}\left( \gsuml{}r_{k-q-1} + \epsilon
     + \sum_{j=1}^{k-q-2}\gamma^{j(\ell m+1)}\left(\gsuml{} r_{k-q-j-1}+\epsilon  \right) \right)\hspace{-1pt}&\{\eqref{eq:1b}\}\\
 &=r_{i}+\gamma\gamma^{(q+1)(\ell m+1)}\left( \sum_{j=0}^{k-q-2}\gamma^{j(\ell m+1)}\left(\gsuml{} r_{k-q-j-1}+\epsilon  \right) \right)\\
&=r_{i}+\gamma\gamma^{q(\ell m+1)}\left( \sum_{j=1}^{k-q-1}\gamma^{j(\ell m+1)}\left(\gsuml{} r_{k-q-j}+\epsilon  \right) \right)\\
&=r_{i}+\gamma\gamma^{q(\ell m+1)}\left(\left(\gsuml{(p+1)} - \gamma^{\ell m}\right)r_{k-q}+ \sum_{j=1}^{k-q-1}\gamma^{j(\ell m+1)}\left(\gsuml{} r_{k-q-j}+\epsilon  \right) \right)\hspace{-15.2pt}\\
&= r_{i}-\gamma^{q(\ell m+1)+1}\gamma^{\ell m}r_{k-q}+\gamma \left(v_k(k+(qm+p+1)\ell) - \indicator{p=0}\gamma^{q(\ell m+1)}\epsilon\right)&\{\eqref{eq:1b}\}\\
&\leq r_{i}-\gamma^{i-k+q}r_{k-q}+\gamma v_k(k+(qm+p+1)\ell)\\
&<  q_{k+1}^\leftarrow(i),
\end{align*}
where we concluded by observing that this is the same result as
\eqref{eq:case1}.
\paragraph{E: In state $i=k+((k-1)m+1)\ell+1$}
\begin{align*}
  q_{k+1}^\leftarrow (i) &=\gamma v_k(i-1) = \gamma v_k(k+((k-1)m+1)\ell)\\
  &=\gamma^{(k-1)(\ell m+1)+1}\left(\frac{\gamma^\ell-\gamma^{\ell(m+1)}}{1-\gamma^\ell}r_1 + \epsilon\right)&\{\eqref{eq:1b}\text{ with $q=k-1$ and $p=0$}\}\\
  &=\gamma^{(k-1)(\ell m+1)+1}\epsilon&\{r_1 = 0\}\\
  &> r_i&\{r_i < 0\}\\
  &= r_i + \gamma v_k(i+\ell-1)& \{v_k(i+\ell+1) = 0 \ \ \eqref{eq:1f}\}\\
  &=q_{k+1}^\rightarrow (i).
\end{align*}
\paragraph{F: In states $i>k+((k-1)m+1)\ell+1$}
Following \eqref{eq:1f} we have $v_k(i-1) = v_k(i+\ell-1) = 0$ and hence
\[
q_{k+1}^\leftarrow (i) = 0 > r_i = q_{k+1}^\rightarrow(i).
\]
\subsubsection{The value function $v_{k+1}$}
In the following we will show that the value function $v_{k+1}$
satisfies Lemma~\ref{lem:policyvalue}. To that end we consider the
value of $((\tpi{k+1,\ell})^m\tpi{k+1} v_k)(s_0)$ by analysing the
trajectories obtained by first following $m$ times $\pi_{k,\ell}$ then
$\pi_{k+1}$ from various starting states $s_0$.
\par
Given a starting state $s_0$ and a non stationary policy
$\pi_{k+1,\ell}$, we will represent the trajectories as a sequence of
triples $(s_i, a_i, r(s_i, a_i))_{i=0,\dots,\ell m}$ arranged in a
``trajectory matrix'' of $\ell$ columns and $m$ rows. Each column
corresponds to one of the policies
$\pi_{k+1},\pi_k,\dots,\pi_{k+2-\ell}$. In a column labeled by policy
$\pi_j$ the entries are of the form $(s_i, \pi_j(s_i), r(s_i, \pi_j(s_i))$;
this layout makes clear which stationary policy is used to select the
action in any particular step in the trajectory. Indeed, in column
$\pi_j$, we have $(s_i,\rightarrow,r_j)$ if and only if $s_i=j$, otherwise each entry is of the form $(s_i, \leftarrow,
0)$. Such a matrix accounts for the first $m$ applications of the
operator $\tpi{k+1,\ell}$. One addional row of only one triple $(s_i,
\pi_{k+1}(s_i), r_{\pi_{k+1}}(s_i))$ represents the final application
of $\tpi{k+1}$. After this triple comes the end state of the trajectory
$s_{\ell m+1}$.
\begin{figure}[ht]
  \[
  \begin{tikzpicture}[mymatrixenv]
    \matrix[mymatrix] (foo)  {
      \qquad\pi_{4}\qquad&\pi_{3}&\qquad\pi_{2}\qquad\\
      (10,\leftarrow, 0)  & (9,\leftarrow,0)&(8,\leftarrow,0)\\
      (7,\leftarrow, 0)   & (6,\leftarrow,0)&(5,\leftarrow,0)\\
      (4,\rightarrow, r_4)& (6,\leftarrow,0)&(5,\leftarrow,0)\\
      (4,\rightarrow, r_4)& (6,\leftarrow,0)&(5,\leftarrow,0)\\
      (4,\rightarrow, r_4)& \boxed{6}&\\
    };
    \mymatrixbraceright[foo]{2}{5}{$m=4$ times}
    \mymatrixbracetop[foo]{1}{3}{$\ell=3$ steps}
  \end{tikzpicture}
  \]
  \caption{The trajectory matrix of policy $\pi_{4,\ell}$ starting from state $10$ with $m=4$ and $\ell=3$.}
  \label{fig:ex-traj-mat}
\end{figure}
\begin{example}
  Figure~\ref{fig:ex-traj-mat} depicts the trajectory matrix of policy
  $\pi_{4,\ell}=\pi_4\pi_3\pi_2$ with $m=4$ and $\ell=3$. The trajectory
  starts from state $s_0=10$ and ends in state $s_{\ell m+1} = 6$. The
  $\leftarrow$ action is always taken with reward $0$ except when in
  state $4$ under the policy $\pi_4$. From this matrix we can deduce
  that, for any value function $v$: 
  \begin{align*}
  ((\tpi{4,\ell})^m\tpi{4} v)(10) &= \gamma^{6}r_4 + \gamma^{9}r_4 + \gamma^{12}r_4 + \gamma^{13}v(6)\\
   &= \gamma^{2\ell}r_4 + \gamma^{3\ell}r_4 + \gamma^{4\ell}r_4 + \gamma^{4\ell+1}v(6)\\
  &= \frac{\gamma^{2\ell} - \gamma^{(m+1)\ell}}{1-\gamma^\ell}r_4 + \gamma^{\ell m+1}v(6).
  \end{align*}
\end{example}
\par
With this in hand, we are going to prove each case of Lemma~\ref{lem:policyvalue} for $v_{k+1}$.
\paragraph{In states $i< k+1$}
Following $m$ times $\pi_{k+1,\ell}$ and then $\pi_{k+1}$ starting from
these states consists in taking the $\leftarrow$ action $\ell m+1$ times
to eventually finish either in state $1$ if $i \leq \ell m + 2$ with
value
\[
v_{k+1}(i) = \gamma^{\ell m +1}v_k(1) +\epsilon_{k+1}(i)= -\gamma^{\ell m +1}\gamma^{(k-1)(\ell m +1)}\epsilon = -\gamma^{k(\ell m +1)}\epsilon
\]
or otherwise in state $i-\ell m - 1 < k$ with value
\[
v_{k+1}(i) = \gamma^{\ell m +1}v_k(i-\ell m-1) +\epsilon_{k+1}(i)= -\gamma^{\ell m
  +1}\gamma^{(k-1)(\ell m +1)}\epsilon = -\gamma^{k(\ell m +1)}\epsilon
\]
This matches Equation~\eqref{eq:1a} in both cases.
\paragraph{In states $i=k+1+(qm+p+1)\ell$}
Consider the states $i=k+1+(qm+p+1)\ell$ with $q\geq 0$ and $0\leq p<m$.
Following $m$ times $\pi_{k+1,\ell}$ and then $\pi_{k+1}$ starting from
state $i$ gives the following trajectories:
\begin{itemize}
\item when $q=0$, (\textit{i.e.} $i=k+1+(p+1)\ell$):
\end{itemize}
\[
\begin{tikzpicture}[mymatrixenv]
    \matrix[mymatrix] (foo)  {
      \qquad\qquad\pi_{k+1}\qquad\qquad&\pi_{k}&\dots&\qquad\qquad\pi_{k-\ell+2}\qquad\qquad\\
      (k+1+(p+1)\ell,\leftarrow,0)&(k+(p+1)\ell,\leftarrow,0)&\dots&(k+p\ell+2,\leftarrow,0)\\
      (k+1+p\ell,\leftarrow,0)&(k+p\ell,\leftarrow,0)&\dots&(k+(p-1)\ell+2,\leftarrow,0)\\
      \vdots&\vdots&\vdots&\vdots\\
      (k+1+\ell,\leftarrow,0)&(k+\ell,\leftarrow,0)&\dots&(k+2,\leftarrow,0)\\
      (k+1,\rightarrow,r_{k+1})& (k+\ell,\leftarrow,0)&\dots&(k+2,\leftarrow,0)\\
      \vdots&\vdots&\vdots&\vdots\\
      (k+1,\rightarrow,r_{k+1})& (k+\ell,\leftarrow,0)&\dots&(k+2,\leftarrow,0)\\
      (k+1,\rightarrow,r_{k+1})&\boxed{k+\ell}&&\\
    };
    \mymatrixbraceright[foo]{2}{5}{$p+1$ times}
    \mymatrixbraceright[foo]{6}{8}{$m-p-1$ times}
    \mymatrixbracetop[foo]{1}{4}{$\ell$ steps}
\end{tikzpicture}
\]
Using \eqref{eq:1b} with $q=p=0$ as our induction hypothesis, this gives
\begin{align*}
\lefteqn{\left((\tpi{k+1,\ell})^m\tpi{k+1} v_k\right)(i) = \sum_{j=p+1}^{m}\gamma^{\ell j} r_{k+1} + \gamma^{\ell m +1}v_k(k+\ell)}\\
&= \sum_{j=p+1}^{m}\gamma^{\ell j} r_{k+1}+ \gamma^{\ell m +1} \left(\gsuml{}r_{k} +
    \epsilon + \sum_{j=1}^{k-1}\gamma^{j(\ell
        m+1)}\left(\gsuml{} r_{k-j}+\epsilon\right) \right)\\
&=\gsuml{(p+1)}r_{k+1} +\sum_{j=1}^{k}\gamma^{j(\ell
        m+1)}\left(\gsuml{} r_{k-j}+\epsilon\right)\\
\end{align*}
Accounting for the error term and the fact that $i=k+1+\ell \iff p=q=0$, we get
\begin{align*}
v_{k+1}(i) &= \left((\tpi{k+1,\ell})^m\tpi{k+1} v_k\right)(i) + \indicator{i=k+1+\ell}\epsilon\\
&=\gsuml{(p+1)}r_{k+1} +\indicator{p=0}\epsilon+\sum_{j=1}^{k}\gamma^{j(\ell
        m+1)}\left(\gsuml{} r_{k-j}+\epsilon\right)\\
\end{align*}
which is \eqref{eq:1b} for $k+1$ and $q=0$ as desired.
\begin{itemize}
\item when $1\leq q \leq k$:
\end{itemize}
In this case we have $i - (\ell m +1) \geq k+1$, meaning that $k+1$, the
first state where the $\rightarrow$ action would be available
is unreachable (in the sense that the
tractory could end in $k+1$, but no action will be taken
there). Consequently the $\leftarrow$ action is taken $\ell m+1$ times
and the system ends in state $i - \ell m -1 =
k+((q-1)m+p+1)\ell$. Therefore, using \eqref{eq:1b} as induction hypothesis and the fact that $i\not\in\{k+1, k+\ell+1
\}\implies \epsilon_{k+1}(i)=0$, we have:
\begin{align*}
    v_{k+1}(i) &= \gamma^{\ell m+1}v_k(k+((q-1)m+p+1)\ell) + \epsilon_{k+1}(i)
  \\
  &=\gamma^{q(\ell m+1)}\left( \gsuml{(p+1)}r_{k+1-q} +
    \indicator{p=0}\epsilon + \sum_{i=1}^{k-q}\gamma^{i(\ell
      m+1)}\left(\gsuml{}r_{k+1-q-k}+\epsilon \right) \right),
\end{align*}
which statisfies \eqref{eq:1b} for $k+1$.
\paragraph{In state $k+1$}
Following $m$ times $\pi_{k+1,\ell}$ and then $\pi_{k+1}$ starting from
$k+1$ gives the following trajectory:

\[
\begin{tikzpicture}[mymatrixenv]
    \matrix[mymatrix] (foo)  {
       \qquad\pi_{k+1}\qquad&\pi_{k}&\dots&\qquad\pi_{k-\ell+2}\qquad\\
      (k+1,\rightarrow,r_{k+1})& (k+\ell,\leftarrow,0)&\dots&(k+2,\leftarrow,0)\\
      \vdots&\vdots&\vdots&\vdots\\
      (k+1,\rightarrow,r_{k+1})& (k+\ell,\leftarrow,0)&\dots&(k+2,\leftarrow,0)\\
      (k+1,\rightarrow,r_{k+1})&\boxed{k+\ell}&&\\
    };
    \mymatrixbraceright[foo]{2}{4}{$m$ times}
    \mymatrixbracetop[foo]{1}{4}{$\ell$ steps}
\end{tikzpicture}
\]

As a consequence, with \eqref{eq:1c} as induction hypothesis we have:
\begin{align*}
\lefteqn{\left((\tpi{k+1,\ell})^m\tpi{k+1} v_k\right)(k+1) =
\frac{1-\gamma^{\ell(m+1)}}{1-\gamma^\ell}r_{k+1}+\gamma^{\ell m
  +1}v_k(k+\ell)}\\
&=
r_{k+1} + \gsuml{} r_{k+1}+ \gamma^{\ell m+1} \left(\gsuml{} r_k + \epsilon + \sum_{j=1}^{k-1}\gamma^{j(\ell m+1)}\left(\gsuml{} r_{k-j} + \epsilon\right)\right)\\
&=
r_{k+1} + \gsuml{} r_{k+1}+\sum_{j=1}^{k}\gamma^{j(\ell m+1)}\left(\gsuml{} r_{k-j+1} + \epsilon\right)\\
&= r_{k+1} + v_{k+1}(k+\ell+1) - \epsilon
\end{align*}
Hence,
\begin{align*}
v_{k+1}(k+1) &= \left((\tpi{k+1,\ell})^m\tpi{k+1} v_k\right)(k+1) + \epsilon_{k+1}(k+1) \\
&=v_{k+1}(k+\ell+1) + r_{k+1} - 2\epsilon,
\end{align*}
which matches \eqref{eq:1c}.

\paragraph{In states $i=k+1+q\ell+p$}
For states $i=k+1+q\ell+p$ with $0\leq q \leq km-1$ and $1\leq p<\ell$,
the policy $\pi_{k+1,\ell}$ always takes the $\leftarrow$ action with either one of the following trajectories
\begin{itemize}
\item when $q\geq m$:
\end{itemize}

\[
\begin{tikzpicture}[mymatrixenv]
    \matrix[mymatrix] (foo)  {
       \qquad\qquad\pi_{k+1} \qquad\qquad&\pi_{k}&\dots& \qquad\qquad\pi_{k-\ell+2} \qquad\qquad\\
      (k+1+q\ell+p, \leftarrow, 0)& (k+q\ell+p, \leftarrow, 0)& \dots&(k+(q-1)\ell+p+2, \leftarrow, 0)\\
      \vdots&\vdots&\vdots&\vdots\\
      (k+1+(q-m+1)\ell+p, \leftarrow, 0)& (k+q\ell+p, \leftarrow, 0)& \dots&(k+(q-m)\ell+p+2, \leftarrow, 0)\\
      (k+1+(q-m)\ell+p, \leftarrow, 0)&\boxed{k+(q-m)\ell+p}&&\\
    };
    \mymatrixbraceright[foo]{2}{4}{$m$ times}
    \mymatrixbracetop[foo]{1}{4}{$\ell$ steps}
\end{tikzpicture}
\]

As a consequence, with \eqref{eq:1d} as induction hypothesis we have:
\[
v_{k+1}(i) = \left((\tpi{k+1,\ell})^m\tpi{k+1} v_k\right)(i) = \gamma^{\ell m+1}v_{k}(k+(q-m)\ell+p) = -\gamma^{\ell m+1}\gamma^{(k-1)(\ell m+1)}\epsilon
= -\gamma^{k(\ell m+1)}\epsilon
\]
which satisfies \eqref{eq:1d} in this case.
\begin{itemize}
\item when $q < m$:
\end{itemize}

Assuming that negative states correspond to state $1$, where the action
is irrelevant, we have the following trajectory:

\[
\begin{tikzpicture}[mymatrixenv]
  \matrix[mymatrix] (foo2)  {
    \qquad\qquad\pi_{k+1}\qquad\qquad&\dots&\qquad\qquad\pi_{k-\ell+2}\qquad\qquad\\
  (k+1+q\ell+p, \leftarrow, 0)& \dots&(k+(q-1)\ell+p+2, \leftarrow, 0)\\
  \vdots&\vdots&\vdots\\
  (k+1+\ell+p, \leftarrow, 0)& \dots&(k+p+2, \leftarrow, 0)\\
  (k+1+p, \leftarrow, 0)& \dots& ( k-\ell+p+2, \leftarrow, 0)\\
  ( k+1-\ell+p, \leftarrow, 0)&\dots&( k-2\ell+p+2, \leftarrow, 0)\\
   \vdots&\vdots&\vdots\\
   ( k+1-(m-q-1)\ell+p, \leftarrow, 0)&\dots&( k-(m-q)\ell+p+2, \leftarrow, 0)\\
   ( k+1-(m-q)\ell+p, \leftarrow, 0)& & \!\!\!\!\!\!\!\!\!\!\!\!\boxed{k+(q-m)\ell+p}\!\!\!\!\\
   };
   \mymatrixbraceright[foo2]{2}{4}{$q$ times}
    \mymatrixbraceright[foo2]{5}{8}{$m-q$ times}
   \mymatrixbracetop[foo2]{1}{3}{$\ell$ steps}
\end{tikzpicture}
\]

In the above trajectory, one can see that only the $\leftarrow$ action is
taken (ignoring state $1$). Indeed, since we follow the
policies $\pi_{k+1}\pi_k,\dots,\pi_{k-\ell+2}$ the $\rightarrow$ action
may only be taken in states $k+1,k,\dots, k-\ell+2$. When state $k+1$ is
reached, the selected action is $\pi_{k-p+1}(k+1)$ which is $\leftarrow$
since $p \geq 1$. The same reasonning applies in the next states
$k,\dots, k-\ell+1$, where $p\geq 1$ prevents to use a policy that would
select the $\rightarrow$ action in those states.


Since $p-\ell < 0$ the trajectory always terminates in a state $j< k$
with value $v_{k}(j) = -\gamma^{(k-1)(\ell m-1)}\epsilon$ as for the $q\geq
m$ case, which allows to conclude that \eqref{eq:1d} also holds in
this case.

\paragraph{In states $i=k+1+km\ell+p$}
Observe that
following $m$ times $\pi_{k+1,\ell}$ and then  $\pi_{k+1}$ once amounts
to always take $\leftarrow$ actions. Thus, one eventually finishes in state
$k+(k-1)m\ell+p\geq k+1$, which, since $\epsilon_k(i) = 0$, gives
\[
v_{k+1}(i) = \left((\tpi{k+1,\ell})^m\tpi{k+1} v_k\right)(i) = \gamma^{\ell m+1}v_{k}(k+(k-1)m\ell+p) = -\gamma^{\ell m+1}0
= 0,
\]
satisfiying \eqref{eq:1e}.

\paragraph{In states $i>k+1+(km+1)\ell$}
In these states, the action $\leftarrow$ is taken $\ell m +1$ times
ending up in state $j>k+((k-1)m+1)\ell$, with value $v_k(j) = 0$, from
which $v_{k+1}(i) = 0$ follows as required by \eqref{eq:1f}.

\section{Empirical Illustration}

\begin{figure}
\includegraphics[width=\textwidth]{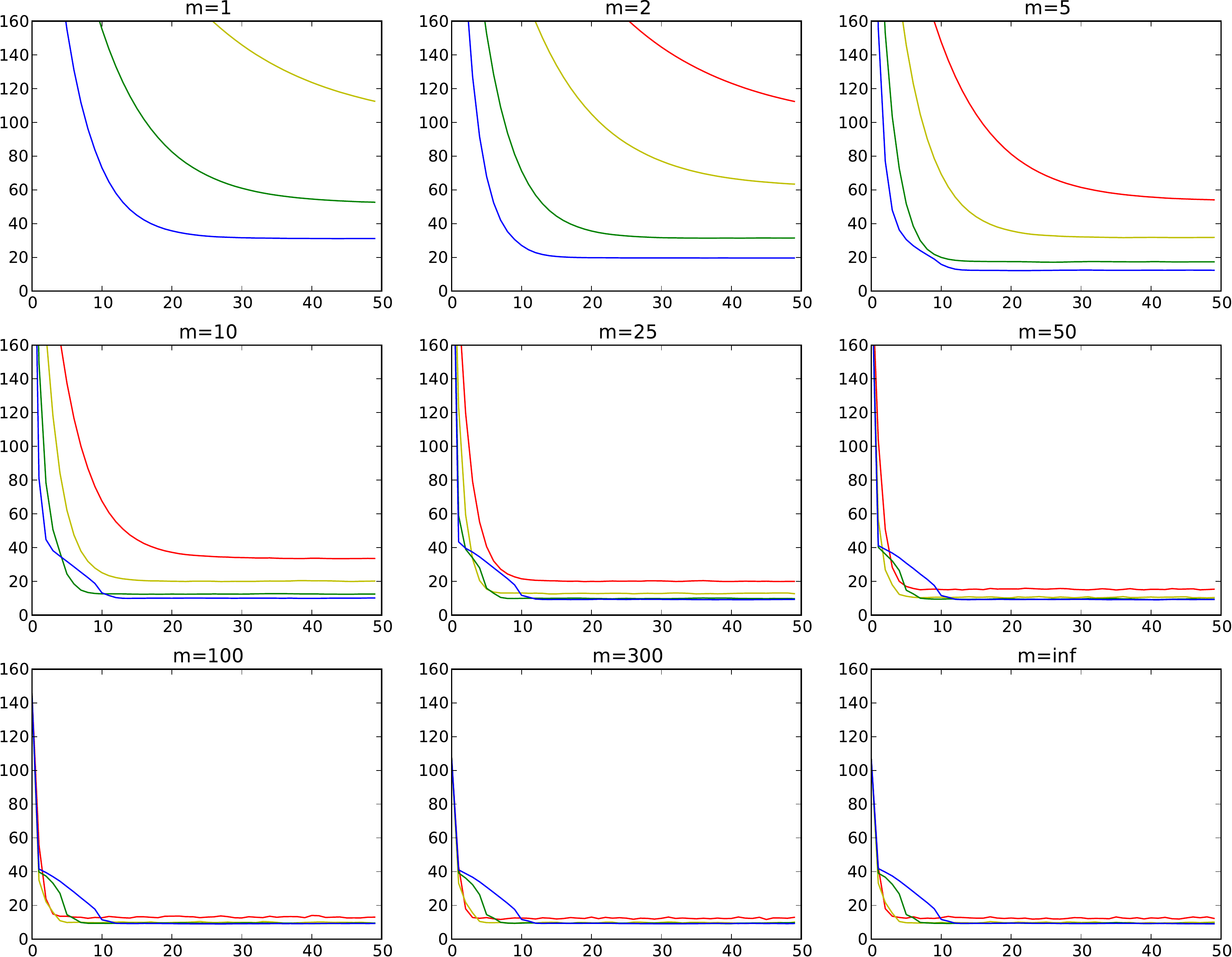}
\caption{Average error of policy $\pi_{k,\ell}$ per iteration $k$ of
  MPI. Red lines for $\ell=1$, yellow for $\ell=2$, green for
  $\ell=5$ and blue for $\ell=10$.}
\label{fig:xp}
\end{figure}

\label{sec:empirical}
In this last section, we describe an empirical illustration of our new
variation of MPI on the dynamic location problem
from~\citet{bertsekas2012q}. The problem involves a repairman moving
between $n$ sites according to some transition probabilities. As to
allow him do his work, a trailer containing supplies for the repair
jobs can be relocated to any of the sites at each decision epoch. The
problem consists in finding a relocation policy for the trailer
according the repairman's and trailer's positions which maximizes the
discounted expectation of a reward function.  \par Given $n$ sites,
the state space has $n^2$ states comprising the locations of both the
repairman and the trailer. There are $n$ actions, each one corresponds
to a possible destination of the trailer. Given an action
$a=1,\dots,n$, and a state $s=(s_r, s_t)$, where the repairman and the
trailer are at locations $s_r$ and $s_t$, respectively, we define the
reward as $r(s,a) = -\abs{s_r - s_t} - \abs{s_t - a}/2$. At any
time-step the repairman moves from its location $s_r < n$ with uniform
probability to any location $s_r\leq s_r'\leq n$; when $s_r=n$, he
moves to site $1$ with probability $0.75$ or otherwise stays. Since
the trailer moves are deterministic, the transition function is
\[
T((s_r,s_t), a, (s_r',a)) = \left\{
\begin{array}{ll}
  \frac 1 {n-s_r+1}&\text{if } s_r < n\\
  0.75 &\text{if } s_r = n \land s_r' = 1\\
  0.25 &\text{if } s_r = n \land s_r' = n\\
\end{array}
\right.
\]
and $0$ everywhere else.
\par

We evaluated the empirical performance gain of using non-stationary policies by implementing the algorithm using random error vectors
$\epsilon_k$, with each component being uniformly random between $0$ and some
user-supplied value $\epsilon$. The adjustable size (with $n$) of the
state and actions spaces allowed to compute an optimal policy to compare
with the approximate ones generated by MPI for all combinations of
parameters $\ell\in\{1, 2, 5, 10\}$ and $m\in\{1, 2, 5, 10, 25,
\infty\}$. Recall that the cases $m=1$ and $m=\infty$ correspond respectively to
the non-stationary variants of VI and PI of \citet{ScherrerLesner2012}, while
the case $\ell=1$ corresponds to the standard MPI algorithm. 
We used $n=8$ locations, $\gamma=0.98$ and $\epsilon=4$ in all experiments.
\par

Figure~\ref{fig:xp} shows the average value of the error
$v_*-v_{\pi_{k,\ell}}$ per iteration for the different values of
parameters $m$ and $\ell$. For each parameter combination, the results
are obtained by averaging over 250 runs. While higher values of $\ell$
impacts computational efficiency (by a factor $O(\ell)$) it always
results with better performance. Especially with the lower values of
$m$, a higher $\ell$ allows for faster convergence. While increasing $m$,
this trend fades to be finally reversed in favor of faster convergence
for small $\ell$. However, while small $\ell$ converges faster, it is
with greater error than with higher $\ell$ after convergence. It can be
seen that convergence is attained shortly after the $\ell^\mathit{th}$
iteration which can be explained by the fact that the first policies (involving $\pi_0, \pi_{-1}, \dots, \pi_{-\ell+2}$), are of poor quality and the algorithm must perform
at least $\ell$ iterations to ``push them out'' of $\pi_{k,\ell}$. 
\par

\begin{figure}
\begin{center}
\includegraphics[width=.9\textwidth]{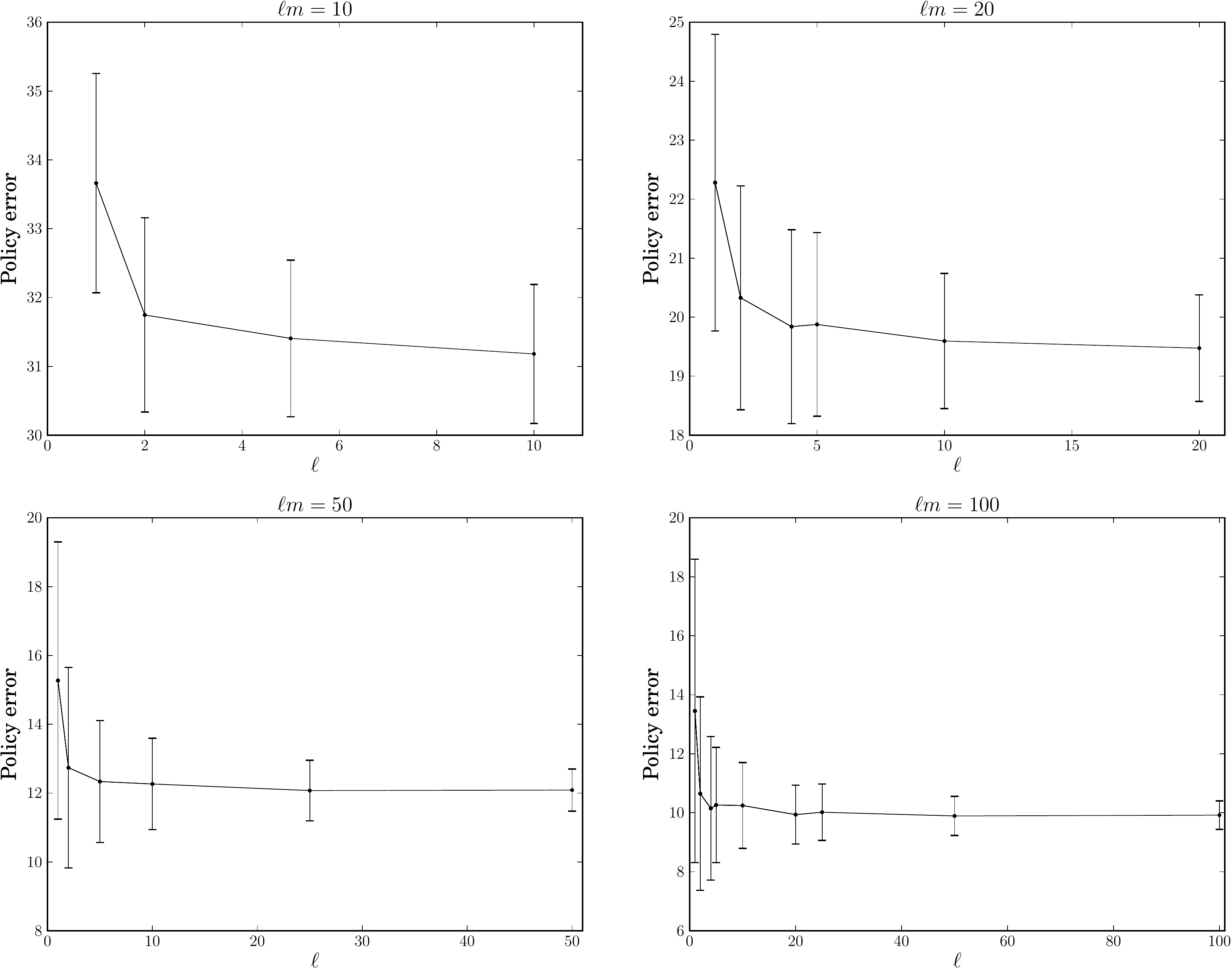}
\end{center}
\caption{Policy error and standard deviation after $150$ iterations for
  different different values of $\ell$. Each plot represents a fixed
  value of the product $\ell m$. Data is collected over $250$ runs with
  $n=8$.}
\label{fig:constant-lm}
\end{figure}

We conducted a second experiment to study the relative influence of the
parameters $\ell$ and $m$. From the observation that the time complexity
of an iteration of MPI can be roughly summarized by the number $\ell
m+1$ of applications of a stationary policy's Bellman operator, we ran
the algorithm for fixed values of the product $\ell m$ and measured the
policy error for varying values of $\ell$ after $150$ iterations. These
results are depicted on Figure~\ref{fig:constant-lm}. This setting
gives insight on how to
set both parameters for a given ``time budget'' $\ell m$.  While runs
with a lower $\ell$ are slightly faster to converge, higher values
always give the best policies after a sufficient number of
iterations. It appears that favoring $\ell$ instead of $m$ seems to
always be a good approach since it also greatly reduces the variance
across all runs, showing that non-stationarity adds robustness to
the approximation noise.



\bibliography{biblio.bib} 

\end{document}